\documentclass[letterpaper, 11pt, oneside]{book}
\usepackage{amssymb,amsmath,amsthm}
\usepackage{latexsym}
\usepackage{geometry}
\usepackage{textcomp}
\usepackage{notoccite}
\usepackage{enumitem}
\usepackage{algorithm}
\usepackage{booktabs}
\usepackage[linktocpage=true,hyperfootnotes=false]{hyperref}
\usepackage{textcomp}
\usepackage[titletoc]{appendix}
\usepackage{setspace} 

\usepackage{amsfonts}
\usepackage[mathscr]{eucal}
\usepackage[dvips]{epsfig,graphics}

\newtheoremstyle{assumptionstyle}{0in}{0in}{\normalfont}{0.5em}{\itshape}{:}{.5em}{}
\newtheoremstyle{propertystyle}{0in}{0in}{\normalfont}{}{\bf}{:}{.5em}{}
\theoremstyle{plain}

\theoremstyle{propertystyle}

\theoremstyle{propertystyle}

\usepackage{algorithm}
\usepackage{algpseudocode}

\usepackage{amsthm}
\usepackage{latexsym}
\usepackage{amsfonts}
\usepackage{amssymb}
\usepackage[mathscr]{eucal}
\usepackage{amsmath}
\usepackage[dvips]{epsfig,graphics}
\usepackage{array}
\usepackage{graphicx}
\usepackage{color}
\usepackage{varwidth}
\usepackage{setspace} 
\usepackage{xspace}
\usepackage{textcomp} 
\usepackage{tabularx}
\usepackage{lscape} 
\usepackage{rotating} 
\usepackage{verbatim}
\usepackage{fancyvrb}
\usepackage{gensymb}
\usepackage{wasysym}

\usepackage{xspace}

 \makeatletter
    \def\thebibliography#1{\chapter*{References\@mkboth
      {REFERENCES}{REFERENCES}}\list
      {[\arabic{enumi}]}{\settowidth\labelwidth{[#1]}\leftmargin\labelwidth
	\advance\leftmargin\labelsep
	\usecounter{enumi}}
	\def\newblock{\hskip .11em plus .33em minus .07em}
	\sloppy\clubpenalty4000\widowpenalty4000
	\sfcode`\.=1000\relax}
    \makeatother

\makeatletter
\def\@makechapterhead#1{%
  \vspace*{50\p@}%
  {\raggedright \normalfont
    \interlinepenalty\@M
    \hspace{-0.375in} \huge\bfseries  \quad #1\par\nobreak
    \vskip 40\p@
  }}
\makeatother

\makeatletter
\renewcommand*\l@chapter[2]{%
  \ifnum \c@tocdepth >\m@ne
    \addpenalty{-\@highpenalty}%
    \vskip 1.0em \@plus\p@
    \setlength\@tempdima{1.5em}%
    \begingroup
      \parindent \z@ \rightskip \@pnumwidth   
      \parfillskip -\@pnumwidth
      \leavevmode \bfseries
      \advance\leftskip\@tempdima
      \hskip -\leftskip
      #1\nobreak\ 
       \leaders\hbox{$\m@th
        \mkern \@dotsep mu\hbox{.}\mkern \@dotsep
        mu$}\hfil\nobreak\hb@xt@\@pnumwidth{\hss #2}\par
      \penalty\@highpenalty
    \endgroup
  \fi}
\makeatother

\begin{document}

\frontmatter
\begin{center}
\thispagestyle{empty}
\vspace*{0.65in}
\begin{singlespace}
{\bf {\Huge{Spatio-Temporal Wind Modeling for UAV Simulations}}}\\
\vspace{0.6 in}
\Large{by Kenan Cole and Adam Wickenheiser}\\
\vspace{0.5in}
\large{\today}

\vspace{4in}
This work originally appeared as a chapter in the dissertation entitled Reactive Trajectory Generation and Formation Control for Groups of UAVs in Windy Environments \cite{ColePhD2018}.
\end{singlespace}
\end{center}

\chapter*{}
\begin{center}
\vspace{2in}
\begin{singlespace}
\begin{table*}[h!]
\begin{tabular}{m{3.25in} m{0.5in} m{2in}} 
This work is licensed under a Creative Commons ``Attribution 4.0" license & \includegraphics[width=1.75in]{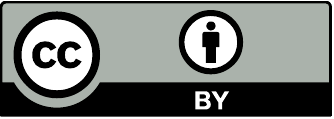} 
\end{tabular}
\end{table*}
\end{singlespace}
\end{center}


\newpage
\addcontentsline{toc}{chapter}{Table of Contents}
\vspace{0.2in}
\begin{singlespace}
\tableofcontents

\clearpage 
\addcontentsline{toc}{chapter}{List of Figures} 
\listoffigures 

\clearpage 
\addcontentsline{toc}{chapter}{List of Tables} 
\listoftables
\end{singlespace} 



\pagebreak
\section*{} 
\addcontentsline{toc}{chapter}{Nomenclature}
\vspace{-0.5in}
\begin{center}
{\bf{Nomenclature}}
\end{center}
\vspace{-0.15in}
\begin{singlespace}
\begin{table*}[h!]
\begin{tabular}{p{0.55in} p{5.2in}}
$A_f$ & experimentally determined parameter used by Forristall model (\ref{SubSecForristall}) \\
$b$ & aircraft wingspan, m (or ft) (\ref{SecOverLandVonKarman}) \\
$B_f$ & experimentally determined parameter used by Forristall model (\ref{SubSecForristall}) \\
$c_p$ & phase speed of the wind waves, m/s (\ref{SecFrictionVelCal}) \\
$C$ & drag coefficient for Harris model (\ref{SubSecHarris}) \\
$D(\omega,\theta)$ & spreading function (\ref{SecSpatialDist}) \\
$D_0$ & coefficient for spreading function (\ref{SecSpatialDist}) \\
$f$ & frequency of wind power spectral density functions, Hz (\ref{SecOverLandVonKarman}) \\
$f^*_f$ & non-dimensional frequency used by Forristall model (\ref{SubSecForristall}) \\
$f^*_h$ & non-dimensional frequency used by Harris model (\ref{SubSecHarris}) \\
$g$ & gravity, 9.81 m/s$^2$ \\
$g_1$ & gust amplitude coefficient (\ref{SubSecGustShape})  \\
$g_{2r}$ & gust time scale coefficient (\ref{SubSecGustShape}) \\
$g_{2f}$ & gust time scale coefficent (\ref{SubSecGustShape}) \\
$g_3$ & gust sigmoid curve error approximation coefficient (\ref{SubSecGustShape}) \\
$g_4$ & pre-gust dip coefficient (\ref{SubSecGustShape})   \\
$g_5$ & post-gust dip coefficient (\ref{SubSecGustShape})\\
$k$ & Von K\'arm\'an constant, 0.42 (\ref{SecTypesOfWind}) \\
$k(\omega)$ & wavenumber (\ref{SecSpatialDist}) \\
$K_{gust,lat}$ & scaling coefficient for lateral spread of gust (\ref{SubSecLatLongDissipation}) \\
$l_x$, $l_y$ & length scale for gust probability of propagation (\ref{SubSecLatLongDissipation})   \\
$l_{gust,span}$ & lateral spread of gust, m (\ref{SubSecLatLongDissipation}) \\
$L$ & scaling length for Harris model, m (\ref{SubSecHarris}) \\
$L_u$ & Von K\'arm\'an model longitudinal length scale, ft (\ref{SecOverLandVonKarman}) \\
$L_v$ & Von K\'arm\'an model lateral length scale, ft (\ref{SecOverLandVonKarman}) \\
$L_w$ & Von K\'arm\'an model vertical length scale, ft (\ref{SecOverLandVonKarman}) \\
$\mathbf{p}$ & position vector in inertial frame \\
$\mathbf{p}_{g}$ & goal position \\
$\mathbf{p}_{gust}$ & gust starting position (\ref{SubSecLatLongDissipation}) \\
$\mathbf{p}^w$ & position vector in wind frame (\ref{SecSpatialDist}) \\
$s$ & spreading parameter (\ref{SecSpatialDist}) \\
$S(\omega)$ & power spectral density function, m$^2$/s (\ref{SubSecPSDCompare}) \\
$t$ & simulation time, s \\
$t_{g1}$ & time span from start of gust to peak, s (\ref{SubSecGustShape}) \\
$t_{g2}$ & time span from peak to end of gust, s (\ref{SubSecGustShape}) \\
$t_h$ & time span for gust at peak magnitude, m/s (\ref{SubSecGustShape}) \\
$u(z)$ & mean wind speed at altitude $z$ \\
$u_{10}$ & mean wind speed at 10m altitude, m/s (\ref{SecTypesOfWind}) \\
$u_{20}$ & mean wind speed at 20ft altitude, ft/s (\ref{SecOverLandVonKarman}) \\
$u_*$ & friction velocity, m/s (\ref{SecTypesOfWind}) \\
$v_{air}$ & mean wind speed, m/s (\ref{SubSecVectorComponents}) \\
$v_{air,spatial}$ & scalar wind speed at some ($x$,$y$) position and time $t$, m/s (\ref{SecSpatialDist}) \\
$v_{gust}$ & maximum gust magnitude, m/s (\ref{SubSecGustAmp}) 
\end{tabular}
\end{table*}
\begin{table*}
\begin{tabular}{p{0.55in} p{5.2in}}
$v_{gust,nom}$ & gust magnitude, m/s (\ref{SubSecGustShape}) \\
$v_{turb}$ & scalar wind speed at some ($x$,$y$) position and time for any propagation direction and model, m/s (\ref{SecSpatialDist})  \\
$\mathbf{v}_{air}$ & translational wind vector in inertial frame, m/s (\ref{SubSecVectorComponents}) \\
$x_{wa}$ & wave age (\ref{SecFrictionVelCal}) \\
$x$ & position along $\mathbf{x}_I$ axis, m (or ft) \\
$\mathbf{x}_B$ & body frame $x$ axis \\
$\mathbf{x}_I$ & inertial frame $x$ axis \\
$\mathbf{x}_W$ & wind frame $x$ axis (\ref{SecSpatialDist}) \\
$y$ & position along $\mathbf{y}_I$ axis, m (or ft) \\
$\mathbf{y}_B$ & body frame $y$ axis \\
$\mathbf{y}_I$ & inertial frame $y$ axis \\
$\mathbf{y}_W$ & wind frame $y$ axis (\ref{SecSpatialDist}) \\
$z$ & altitude, m (or ft) \\
$\mathbf{z}_B$ & body frame $z$ axis \\
$\mathbf{z}_I$ & inertial frame $z$ axis \\
$\mathbf{z}_W$ & wind frame $z$ axis (\ref{SecSpatialDist}) \\
$z_0$ & sea roughness, m (\ref{SecTypesOfWind}) \\
$z_0^*$ & non-dimensional sea roughness (\ref{SecFrictionVelCal}) \\
$\delta x$, $\delta y$ & distance from vehicle to gust center (\ref{SubSecLatLongDissipation}) \\
$\theta$ & spreading angle, $\pm \frac{\pi}{2}$ (\ref{SecSpatialDist}) \\
$\rho$ & density of water or air, kg/m$^3$ \\
$\sigma$ & wind speed variance for Forristall model, m$^2$/s$^2$ (\ref{SubSecForristall}) \\
$\sigma_{air}$ & standard deviation of wind speed over 10 minute period, m/s (\ref{SubSecGustAmp}) \\
$\sigma_u$ & Von K\'arm\'an model longitudinal intensity, ft/s (\ref{SecOverLandVonKarman}) \\
$\sigma_v$ & Von K\'arm\'an model lateral intensity, ft/s (\ref{SecOverLandVonKarman}) \\
$\sigma_w$ & Von K\'arm\'an model vertical intensity, ft/s (\ref{SecOverLandVonKarman}) \\
$\tau$ & shear stress, MPa (\ref{SecTypesOfWind}) \\
$\tau_{gust}$ & gust timespan, s (\ref{SubSecGustDuration}) \\
$\Phi_f$ & Forristall model power spectral density, m$^2$/s (\ref{SubSecForristall}) \\
$\Phi_h$ & Harris model power spectral density, m$^2$/s (\ref{SubSecHarris}) \\
$\Phi_{os}$ & Ochi and Shin model power spectral density, m$^2$/s (\ref{SubSecOchiShin}) \\
$\Phi_{os}^*$ & Non-dimensional power spectral density for Ochi and Shin model (\ref{SubSecOchiShin})  \\
$\Phi_p$ & Von K\'arm\'an model pitch power spectral density, 1/s (\ref{SecOverLandVonKarman}) \\
$\Phi_q$ & Von K\'arm\'an model roll power spectral density, 1/s (\ref{SecOverLandVonKarman})\\
$\Phi_r$ & Von K\'arm\'an model yaw power spectral density, 1/s (\ref{SecOverLandVonKarman}) \\
$\Phi_u$ & Von K\'arm\'an model longitudinal power spectral density, ft$^2$/s (\ref{SecOverLandVonKarman}) \\
$\Phi_v$ & Von K\'arm\'an model lateral power spectral density, ft$^2$/s (\ref{SecOverLandVonKarman}) \\
$\Phi_w$ & Von K\'arm\'an model vertical power spectral density, ft$^2$/s (\ref{SecOverLandVonKarman}) \\
$\psi$ & randomly generated phase angle from 0 to $2\pi$, rad (\ref{SubSecPSDCompare}) \\
$\omega$ & frequency of wind power spectral density functions, rad (\ref{SecOverLandVonKarman}) \\
$\boldsymbol\omega_{air}$ & rotational wind vector in inertial frame, rad/s (\ref{SubSecVectorComponents})
\end{tabular}
\end{table*}
\end{singlespace}


\mainmatter
\chapter[Introduction]{1 Introduction}
Wind affects the stability and maneuverability of UAVs, which can be particularly dangerous when operating near obstacles or each other. In order to test the effectiveness of formation control laws and the impact of windy environments on the vehicles, spatio-temporal wind fields must be modeled. Each vehicle within the formation experiences unique wind conditions, but these conditions are correlated to the conditions experienced by the other vehicles. 

A realistic wind field is not static; rather it is should include turbulence and gusting. {\it{Turbulence}} is defined by the National Weather Service (NWS) as an ``abrupt or irregular movement of air that creates sharp, quick updrafts or downdrafts" \cite{NWSTurbulence2006}. These perturbations on the mean wind speed are always present with varying but bounded magnitudes and act in any translational or rotational direction. {\it{Gusting}} on the other hand is defined by NWS as ``a sudden, brief increase in speed of the wind" \cite{NWSGust2018}. Gusts are not always present, have a much larger magnitude over the mean wind speed than turbulence, and may propagate through the wind field before dissipating. Both of these phenomena impact the vehicle's ability to follow a trajectory and/or maintain a formation. 

\begin{figure}
	\begin{center}
		\includegraphics[width=3in]{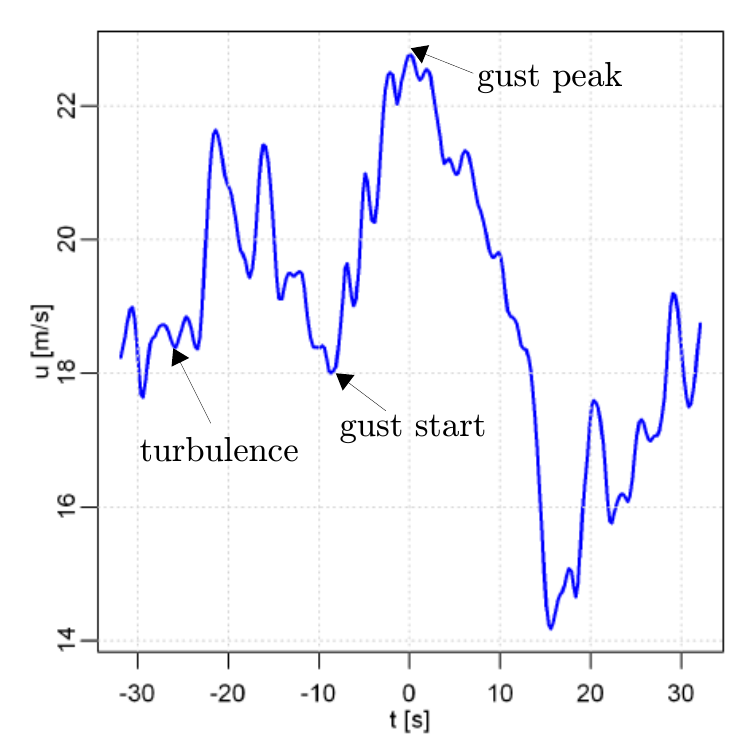}
		\caption[Example wind data from a wind farm in the Netherlands \cite{Branlard2009} (Figure 3.3a).]{\label{FigExIntroRealWind} Example wind data from a wind farm in the Netherlands \cite{Branlard2009} (Figure 3.3a). Notice the turbulence (small perturbations about the mean values) and the gusts (large changes in the wind speed above the mean) in the data. UAVs operating outdoors may be subject to wind disturbances similar to this.}
	\end{center}
\end{figure}

While there are many studies that addresses UAV flight control and formations, there is a much smaller subset that considers the effects of wind on the vehicles. Of the works that do, their wind models are simple compared to real wind data that exhibit turbulence and gusting, such as the wind data from a wind farm in the Netherlands \cite{Branlard2009}, shown in Fig.~\ref{FigExIntroRealWind}. UAV studies such as \cite{Jackson2008}-\nocite{Nelson2006}\nocite{Ratnoo2011}\nocite{Brezoescu2011}\cite{Cabecinhas2013} consider wind as a constant disturbance to their controllers. Similarly, \cite{Bezzo2016} use the maximum velocity with a bounded acceleration as the worst case scenario for their control law evaluation, and \cite{Zhang2015} assume slow velocity changes and therefore model the wind as a piecewise constant function in acceleration. These models do not reflect the turbulence or gusting that is present in a real environment. A slightly more complex wind model of the form $A \sin (\omega t + \phi)$ is used in other studies \cite{McGee2006}\nocite{Escareno2013}-\cite{Tanyer2016}, but like the constant disturbance model it is not stochastic, so it is not a good reflection of realistic turbulence. Moving away from constant and periodic disturbances, \cite{Miller2011,Chen2016} use white noise as a disturbance for UAV simulations. The power spectral density (PSD) of white noise is uniform over all frequencies, which does not reflect measured wind data. The Dryden turbulence model \cite{MILF8785C1980} is utilized in \cite{Waslander2009},\cite{Hancer2010},\cite{Sydney2013},\cite{Kun2016} for testing vehicle control laws; however, this is a stationary process so gusting is not considered. Both turbulence and gusting should be considered in the wind model since they affect the vehicle response and control design. The controller must continuously combat turbulence, whereas it only occasionally experiences gusts, but the gust may exceed the vehicle's control authority. A controller that is robust to one may not be robust to the other, and the vehicle should be capable of handling both. 

In addition to wind generation, the spatial distribution to model the distinct yet correlated conditions experienced by each vehicle is needed. Spatial distributions are considered for a few different applications depending on the objective of the control algorithm. For example, \cite{Lan2016} generates a spatio-temporal wind field using basis functions and time-varying coefficients to test UAVs' abilities to generate trajectories for sampling, estimating, and building a map of an unknown wind field. Their application is focused on building a map over a large area over longer time scales than what is considered in this dissertation. The longer time scales may miss important turbulence and/or potentially significant gust accelerations that are important for controller evaluation. More similar to this research are \cite{Peterson2010},\cite{Peterson2013},\cite{Peterson2015},\cite{Mellish2010},\cite{Sydney2014}, where trajectories are generated to efficiently move individual vehicles or formations of vehicles through wind fields. The flow fields generated are both uniform and nonuniform; however, in all cases no turbulence is considered. Additionally, of the spatial wind field works reviewed none consider gusting.

\begin{figure}
	\begin{center}
		\includegraphics[width=6in]{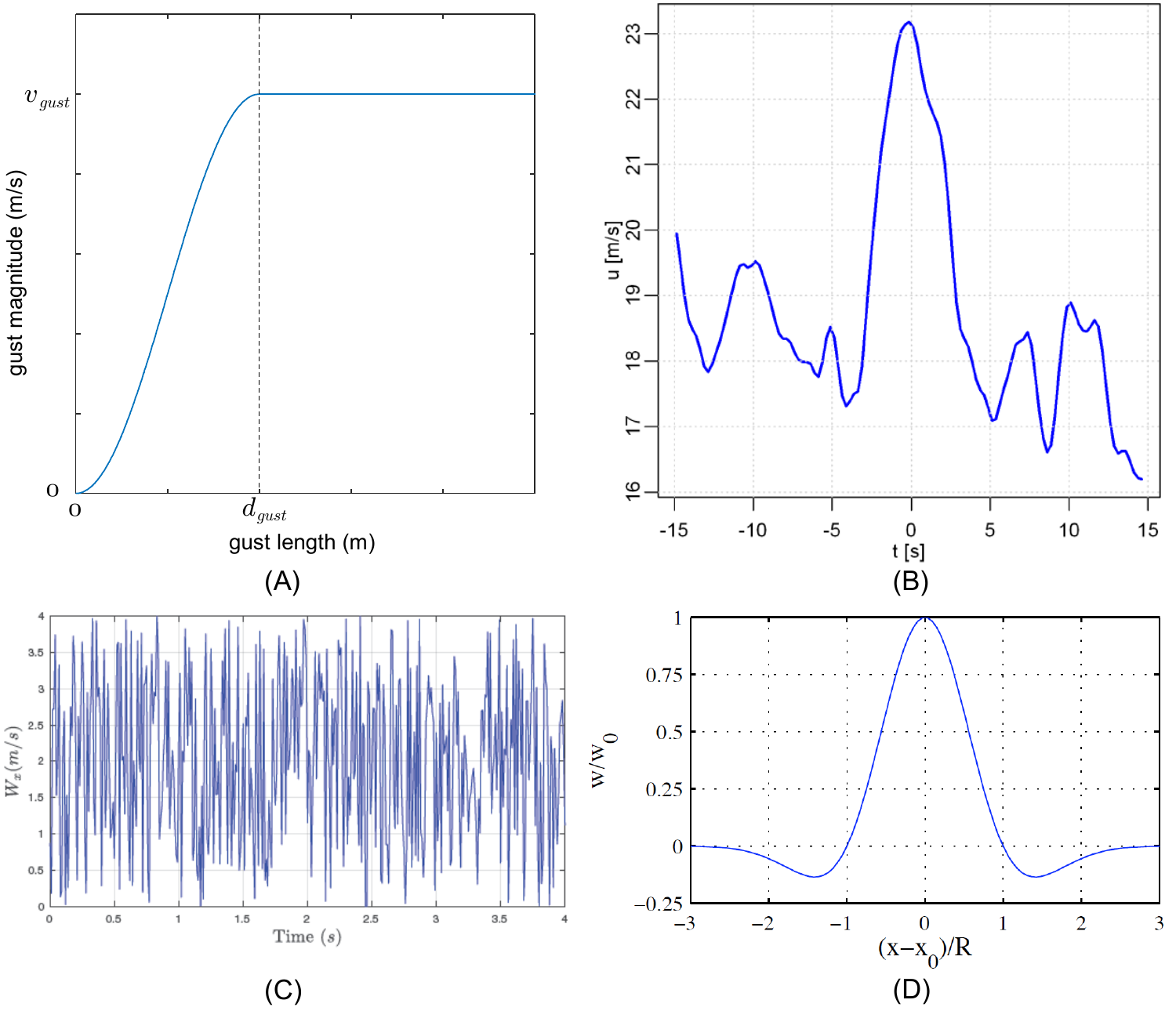}
		\caption[Example gusts from literature]{\label{FigExIntroGusts} Example gusts from the literature. (A) Discrete $(1-\cos)$ gust model defined in \cite{MILF8785C1980}. (B) Experimental gust data from wind farm \cite{Branlard2009} (Figure 3.2). (C) Random data used to represent wind \cite{Fernandez2017} (Figure 5a). (D) Normalized wind profile for testing a soaring controller \cite{Langelaan2008} (Figure 6).}
	\end{center}
\end{figure}

Gust modeling is prominent in wind turbine literature and thermal soaring applications but is less prevalent in UAV control literature. Studies such as \cite{Etele2006} indicate the need for gust and turbulence modeling for UAV applications and provide examples for generating a $(1-\cos)$ curve model that is defined in \cite{MILF8785C1980}. The gust and turbulence are generated for application to a single vehicle with no spatial considerations. Fern\'andez et al.~\cite{Fernandez2017} use random data for the $x$ and $y$ wind velocities, where any increase above the mean is considered a ``gust" by the authors. They use a fan for experimental results but do not compare either the simulated or fan gust models to any real-world data. Alternately, \cite{Langelaan2008} use a gust shape that is similar to the gust shape presented by \cite{Branlard2009} to model vertical gusts for thermal soaring controllers. In both of these applications only a single vehicle is considered, and there are no spatial variations in wind. Examining gusts from a different perspective, \cite{Branlard2009} published a study summarizing the statistics of gusts after recording data over a three year period at a rate of 5 Hz (to capture the high frequency, short time scale variations) at a wind farm in the Netherlands. This work provides useful correlations for gust and wind parameters such as magnitude, duration, and propagation, but does not endorse a single method for spatio-temporal gust creation. Figure \ref{FigExIntroGusts} shows some of the example gust shapes from the literature.

Of the works reviewed, there is no one model that generates a spatio-temporal field that includes turbulence and gusting to use to evaluate UAV controller performance. This chapter addresses this deficiency by developing a model for both over-land and over-water environments. In addition to providing a representative wind field, the model should also be able to run on a personal computer for evaluating a variety of different simulation scenarios. In the following sections, the type of model (i.e.~analytical, large-eddy simulation, stochastic) that best satisfies the requirements for this application is addressed, including the corresponding over-land and over-water turbulence models. Next, the spatio-temporal spreading functions to address how each vehicle experiences unique wind conditions that are correlated with what the other vehicles are experiencing are developed. Lastly, the gusting models are developed and how they interact with the spatio-temporal field is presented. 

\chapter[Types of Wind Models]{2 Types of Wind Models}
\label{SecTypesOfWind}
Environmental wind, turbulence, and gusting models range from very simple analytical models to computational fluid dynamics (CFD) simulations. There are several different types of CFD models, including Reynolds-Averaged Navier-Stokes (RANS), large-eddy simulation (LES), and detached eddy simulation (DES). In all of these CFD models, the area or volume of the fluid is discretized and the equations are solved iteratively. 

\begin{figure}[h!]
	\begin{center}
		\includegraphics[width=6in]{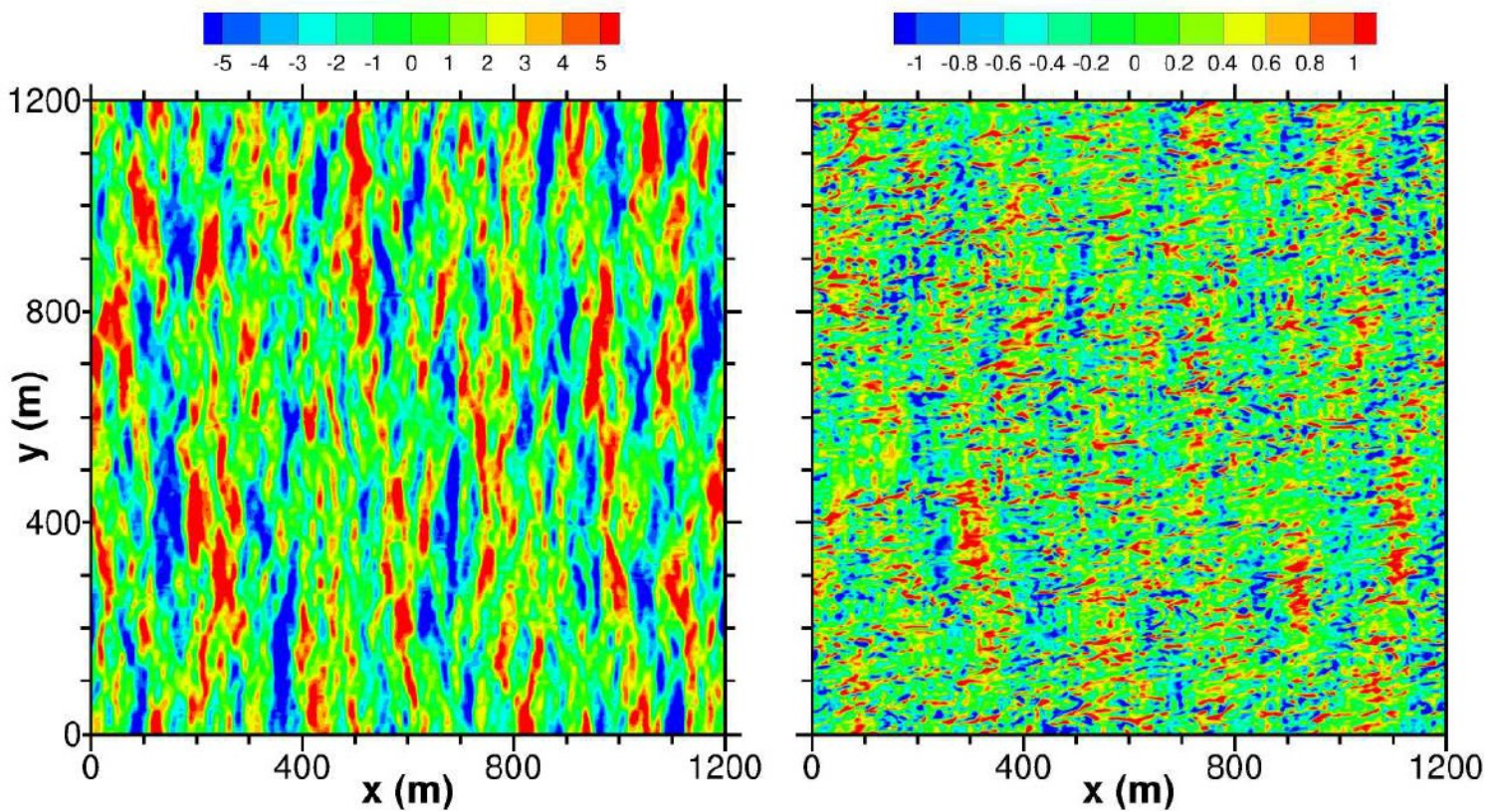}
		\caption[Spatial wind distribution at $z$ = 2.5 m showing the variation in the wind]{\label{figSullivanSpatialWind} Spatial wind distribution at $z$ = 2.5 m showing the variation in the wind over water for heavier wind and waves (left) and much calmer wind and waves (right) \cite{Sullivan}.}
	\end{center}
\end{figure}

\begin{figure}[h!]
	\begin{center}
		\includegraphics[width=6in]{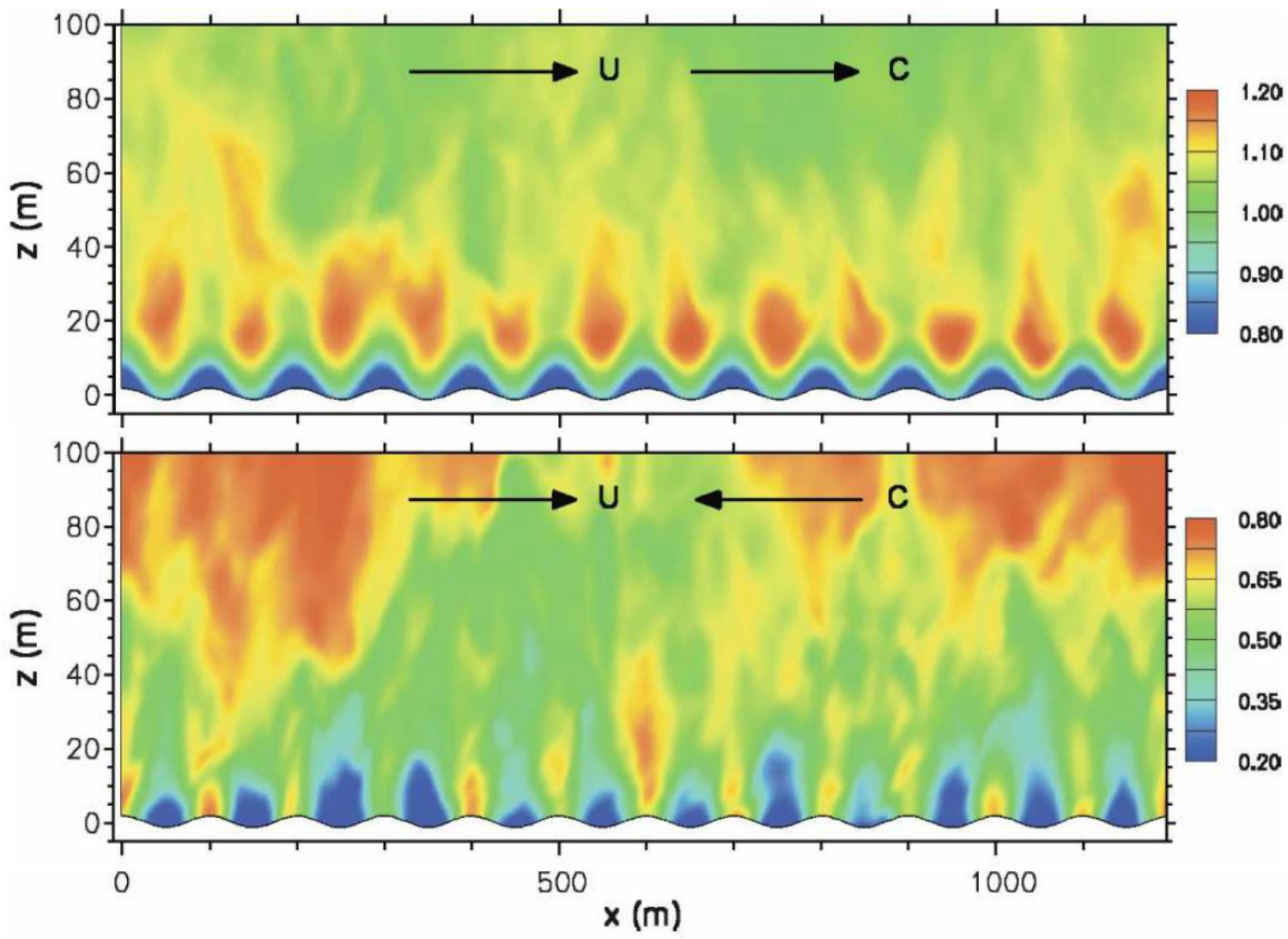}
		\caption[Vertical wind distribution showing the variation when wind ($u$) and waves]{\label{figSullivanVerticalWind} Vertical wind distribution showing the variation when wind ($u$) and waves ($c$) are in the same direction (upper), and opposite directions (lower) \cite{Sullivan2008}.}
	\end{center}
\end{figure}

For the case of wind-over-water modeling, there are several groups that have used large-eddy simulations \cite{Sullivan}\nocite{Sullivan2008}\nocite{Suzuki2010}-\cite{Harcourt2007} to model and understand the wind and wave interactions. Two examples in Figs.~\ref{figSullivanSpatialWind} and \ref{figSullivanVerticalWind} highlight example data from Sullivan, et al.~\cite{Sullivan} where the spatial spread at a constant altitude and the variation of the wind with altitude are shown, respectively. These examples show the capabilities of CFD models to generate realistic and accurate spatio-temporal wind data at high spatial resolutions. Accurate CFD solutions, however, are highly dependent on the quality of the mesh (or grid) as well as the analysis time step. If the problem and initial conditions are not defined appropriately, then the results may be inaccurate and/or may even diverge\cite{Yee1999,Blazek2005}. Additionally, the computational requirements to generate wind data for evaluating formation and vehicle controllers exceeds the computing power of a personal computer, which means CFD is not feasible for use in this work. 

\begin{figure}
	\begin{center}
		\includegraphics[width=3in]{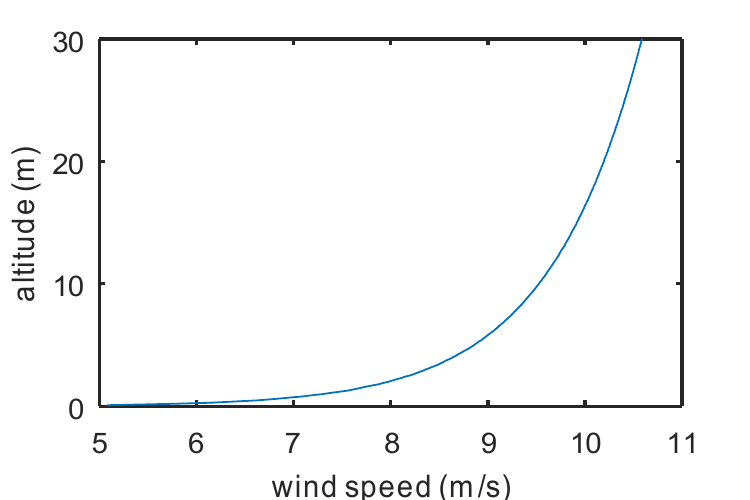}
		\caption[Logarithmic wind profile based for $u_{10}$ = 10 m/s]{\label{figRugglesLogWind} Logarithmic wind profile for $u_{10}$ = 10 m/s (the wind speed at an altitude of 10 m), $u_* = 0.16$ m/s, and $z_0 = 0.0005$ m.}
	\end{center}
\end{figure}

On the opposite end of the spectrum, there are analytical models that describe averaged horizontal wind profiles as a function of altitude \cite{Ruggles1970,Tambke2001}, but offer no time or spatial variations. One such model is a logarithmic profile presented in Fig.~\ref{figRugglesLogWind}, where the dependence with altitude is governed by
\begin{equation}
	\label{EqLogWindProfile}
	u(z) = \frac{u_*}{k} \ln \frac{z}{z_0}
\end{equation}

\noindent where $u_*$ is the friction velocity, $k$ is the Von K\'arm\'an constant taken as 0.42, $z$ is the altitude, and $z_0$ is the sea roughness. The friction velocity is the velocity at the boundary layer between the wind and water and is sometimes called the shear velocity. If shear stress data is available, then the friction velocity is calculated by
\begin{equation}
	\label{EqFrictionVelocity}
	u_* = \sqrt{\frac{\tau}{\rho}}
\end{equation}
\noindent where $\tau$ is the shear stress and $\rho$ is the density. If no shear stress data is available, the friction velocity is solved iteratively with the sea roughness, as discussed in Sec.~\ref{SecFrictionVelCal}.

The desired model lies between the extremes of the analytical and CFD models such that enough of the fidelity is preserved to realistically approximate the wind environment for a group of vehicles, but the computational expense is reduced so that wind can be generated for vehicle formation simulations running on a personal computer. A model that achieves this is a stochastic turbulence model that is expanded into a spatio-temporal distribution by spreading functions (as discussed in Sec.~\ref{SecSpatialDist}). The benefits of a stochastic turbulence model are that the frequency content is preserved and realistic, it is computationally feasible, and it can be combined with a gust model that propagates through the spatio-temporal field. Since it is a stochastic process, the flight controllers cannot be tuned to a specific wind time history. Additionally, since the frequency content is preserved, it can be used to determine the expected response of the vehicle. There are several stochastic turbulence models in the literature for both over-land and over-water environments. These models are examined more closely in Secs.~\ref{SecOverLandVonKarman} to \ref{SubSecOchiShin} to choose one that is appropriate for this research. 

\chapter[Over-Land Wind Turbulence Modeling]{3 Over-Land Wind Turbulence Modeling}
\label{SecOverLandVonKarman}
The over-land case is considered first, where the Dryden and Von K\'arm\'an turbulence models are the standards for modeling wind turbulence. These two models are defined in the military specification MIL-F-8785C \cite{MILF8785C1980} and military handbook MIL-HDBK-1797 \cite{MILHDBK17971997} for aircraft that are flying at both high ($>$ 2000 ft or 610 m) and low altitudes ($<$ 2000 ft or 610 m). While the models are similar in their power spectral densities and the length scale parameters, the Von K\'arm\'an turbulence model better matches data of continuous turbulence \cite{Hoblit1988} and is the preferred model of MIL-F-8785C. 

While the Von K\'arm\'an and Dryden models were developed for large aircraft, Patel et al.~\cite{Patel2008} and Watkins and Vino \cite{Watkins2004,Watkins2006} verified that these models are still reasonable for small UAVs at lower altitudes. Patel et al.~collected data using a small UAV on what they characterize as a ``moderately gusty" day to use for comparison with the Dryden model. They use the change in energy, $dE/dt$, as the measure of a ``gust" and compare the PSD of this data to the Dryden PSD. The slopes and amplitudes match well; the only difference is there are slightly larger amplitudes in the low-frequency range compared to the Dryden model. Similarly, Watkins and Vino \cite{Watkins2004,Watkins2006} recorded wind data with a fixed anemometer on a vehicle in moving and stationary experiments to determine how the turbulence and wind velocities vary for aerial vehicles that are in the low altitude range ($<$ 610 m). They also found good agreement with the low-altitude Von K\'arm\'an model and, similar to Patel et al., saw slightly higher amplitudes at the lower frequencies. The frequency range of interest is 0-8 Hz, where the upper limit is cited by Gage \cite{Gage2003} based on the highest frequency for which the Von K\'arm\'an and Dryden model transfer functions are valid according to \cite{MILHDBK17971997}.

Since there is good agreement for both Von K\'arm\'an and Dryden for smaller UAVs and the military specification states that the Von K\'arm\'an model shall be used when feasible (MIL-F-8785C Section 3.7.1), the Von K\'arm\'an model is chosen for over-land simulations. Table \ref{TblVKLAM} summarizes the low-altitude Von K\'arm\'an model definition and the definitions for the length scales for the low altitude model, since the UAVs are operating below 610 m. The variables are defined in Tbl.~\ref{TblVKVars}, which differentiates what is known and what is calculated to compute the power spectral density.

\begin{table}
	\begin{center}
	\caption[Low-altitude Von K\'arm\'an model]{\label{TblVKLAM}Low altitude Von K\'arm\'an model scaling lengths, intensities, and spectral density}
		\begin{tabular}{|l|c|c|c|} \hline
		& {\bf{Scaling length}} & {\bf{Intensity}} & {\bf{Spectral density }} \\ \hline \hline
		Lon & $L_u = \frac{z}{\left(0.177 + 0.000823 z\right)^{1.2}}$ & $\sigma_u = \frac{\sigma_w}{\left(0.177 + 0.000823z\right)^{0.4}}$ & $\Phi_u = \frac{2 \sigma_u^2 L_u}{\pi u_{20}}\frac{1}{\left(1+\left(1.339 \frac{L_u\omega}{u_{20}}\right)^2\right)^{5/6}}$ \\ \hline
		Lat & $L_v = L_u$ & $\sigma_v = \sigma_u$ & $\Phi_v = \frac{\sigma_v^2 L_v}{\pi u_{20}}\frac{1 + 8/3\left(1.339 \frac{L_v \omega}{u_{20}}\right)^2}{\left(1 + \left(1.339 \frac{L_v \omega}{u_{20}}\right)^2\right)^{11/6}}$ \\ \hline
		Vert & $L_w = z$ & $\sigma_w = 0.1 u_{20}$ & $\Phi_w = \frac{\sigma_w^2 L_w}{\pi u_{20}} \frac{1 + 8/3\left(1.339 \frac{L_w \omega}{u_{20}}\right)^2}{\left(1 + \left(1.339 \frac{L_w \omega}{u_{20}}\right)^2\right)^{11/6}}$ \\ \hline
		Roll & & & $\Phi_p = \frac{\sigma_w^2}{L_w u_{20}}\frac{0.8\left(\frac{\pi L_w}{4b}\right)^{1/3}}{1+\left(\frac{4b\omega}{\pi u_{20}}\right)^{2}}$  \\ \hline
		Pitch & & & $\Phi_q = \frac{\pm\left(\frac{\omega}{u_{20}}\right)^2}{1 + \left(\frac{4b\omega}{\pi u_{20}}\right)^2}\Phi_w(\omega)$ 
		\\ \hline
		Yaw & & & $\Phi_r = \frac{\pm\left(\frac{\omega}{u_{20}}\right)^2}{1 + \left(\frac{3b\omega}{\pi u_{20}}\right)^2}\Phi_v(\omega)$ \\ \hline
		
		\end{tabular}
	\end{center}
\end{table}

\begin{table}
	\begin{center}
		\caption{\label{TblVKVars} Variable descriptions for the Von K\'arm\'an turbulence model}
		\begin{tabular}{|p{0.8in}|p{2.375in}|p{1in}|p{1in}|} \hline
		{\bf{Variable}} & {\bf{Description}} & {\bf{Known/}} & {\bf{Value}} \\ 
		& & {\bf{Calculated}} & \\ \hline \hline
		$\Phi_u$, $\Phi_v$, $\Phi_w$ & PSD for longitudinal, lateral, and vertical turbulence (ft$^2$/s)& Calculated & Table \ref{TblVKLAM} \\ \hline
		$\Phi_p$, $\Phi_q$, $\Phi_r$ & PSD for pitch, roll, and yaw turbulence (1/s) & Calculated & Table \ref{TblVKLAM}\\ \hline
		$\sigma_u$, $\sigma_v$, $\sigma_w$ & Root mean square intensities for longitudinal, lateral, and vertical turbulence velocities (ft/s) & Calculated & Table \ref{TblVKLAM} \\ \hline
		$L_u$, $L_v$, $L_w$ & Scaling lengths for longitudinal, lateral, and vertical turbulence (ft) & Calculated & Table \ref{TblVKLAM} \\ \hline
		$\omega = 2\pi f$ & Frequency (rad) & Known & $f=$0-8 Hz \\ \hline 
		$z$ & UAV altitude (ft) & Known & UAV Dependent \\  \hline
		$u_{20}$ & Mean wind speed at 20 ft (6 m) altitude (ft/s) & Known & Simulation Dependent \\ \hline
		$b$ & Aircraft wingspan (ft) & Known & $b = 1.12$ ft (0.34 m)\\ \hline	
		\end{tabular}
	\end{center}
\end{table}

\begin{figure}
	\begin{center}
		\includegraphics[width=6in]{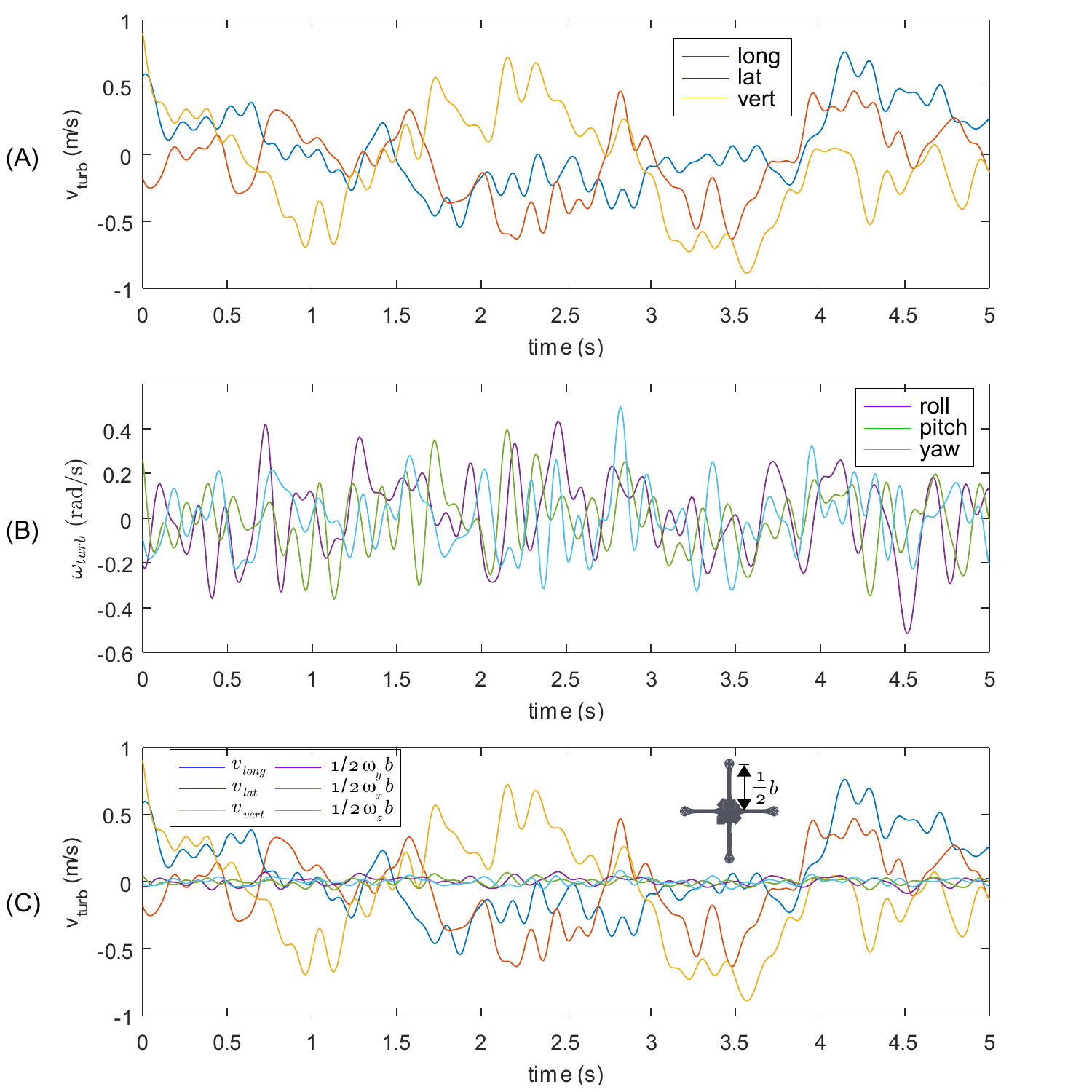}
		\caption[Example time series of the six turbulence components]{\label{FigVKCorrelatedTurbulence} Example time series generation of the six turbulence components generated by the Von K\'arm\'an model with $u_{20}=6$ m/s, $z=10$ m, and $b=0.34$ m. (A) The three translational turbulence components that are added to the mean wind speed along each direction. (B) The rotational components about each axis. (C) Comparison between the three translational components and the translational velocity induced at the motor/propeller assembly from the rotational turbulence components. Notice that the induced translational velocity is comparatively small for a quadrotor since the wingspan is very short.}
	\end{center}
\end{figure}

The Von K\'arm\'an model PSD is used to generate a spatio-temporal wind field as discussed in Sec.~\ref{SecSpatialDist} for use in the over-land UAV simulations. Figure \ref{FigVKCorrelatedTurbulence} shows example time series data for the translational and rotational disturbances. 

\chapter[Over-Water Wind Turbulence Modeling]{4 Over-Water Wind Turbulence Modeling}
The over-water environment is considered next, where there are several stochastic models \cite{Forristall1988}\nocite{Myrhaug2009}\nocite{Brown1991}\nocite{Ochi1988}\nocite{Walmsley1988}\nocite{Islam2009}-\cite{Zaheer2009} that are primarily used in the field of offshore engineering. These models are used to analyze the loads on offshore structures and consider a wide range of environmental conditions including high winds and hurricanes. They are based on and/or verified by data taken from various bodies of water. To choose the best model, popular models such as the Harris \cite{Harris1970}, Forristall \cite{Forristall1988}, and Ochi and Shin \cite{Ochi1988,Fossen1994} models are examined in Secs.~\ref{SubSecHarris} to \ref{SubSecOchiShin} and then compared in Sec.~\ref{SubSecPSDCompare}.

One key assumption in all of the models is a neutrally stable atmosphere. This means that the air is well mixed so any temperature differences do not have an effect. Additionally, neutral stability assumes that the temperature of the water is less than the temperature of the air \cite{Hsu2004}. If this case were reversed, then an unstable atmosphere would result. Walmsley \cite{Walmsley1988} presents a method for incorporating temperature data through use of the Monin-Obukhov length \cite{Obukhov1971}, which requires absolute and potential temperature data. Hsu and Blanchard \cite{Hsu2004} also discuss the impact that an unstable atmosphere has on calculating the profiles through knowledge of the temperature distribution. Since temperature data would be generated somewhat arbitrarily to feed these alternate models, there is minimal (if any) improvement in accuracy to the model, so the assumption of a neutrally stable atmosphere is used.

\section{Harris Model}
\label{SubSecHarris}
Published in 1971, the Harris model was one of the earliest over-water stochastic models proposed \cite{Harris1970}. It is used for high velocity wind, and when compared to hurricane data the Harris spectrum follows the same shape but shows some inaccuracies at particular frequencies. The Harris model PSD is given by
\begin{equation}
	\label{EqHarrisModel}
	\Phi_h = \frac{4 C u_{10}^2 f^*_h}{f \left(2+f^{*^2}_h\right)^{5/6}}
\end{equation}
\noindent where $C$ is the drag coefficient, $f^*_h$ is the non-dimensional frequency, and $u_{10}$ is the wind speed at 10 m altitude. The drag coefficient is a constant equal to 0.002 in what is described qualitatively as ``rough" wave conditions and 0.0015 in ``moderate" conditions \cite{Forristall1988}. The non-dimensional frequency for the Harris model is defined by
\begin{equation}
	\label{EqHarrisNonDimFreq}
	f^*_h = \frac{f L}{u_{10}}
\end{equation}
\noindent where $f$ is the frequency in Hz and $L$ is a scaling length. The scaling length was initially defined as 1800 m to fit the available data \cite{Forristall1988,Harris1970}. Harris \cite{Harris1970} later corrected this value by determining that the length scales increase with altitude and defined a length scale range from approximately 50 m to 126 m for altitudes from 18 m to 182 m. This range of length scales is closer to what Chakrabarti \cite{Chakrabarti1990} presents in his summary of Harris's model in which he defines scaling lengths ranging from 60 m to 400 m. 

Interestingly, the model in Eq.~\ref{EqHarrisModel} was an improvement of another spectral model, the Davenport model \cite{Davenport1961}, and has the same form as the Von K\'arm\'an model summarized in Tbl.~\ref{TblVKLAM}. The coefficients differ as Harris matched experimental wind data and Von K\'arm\'an matched wind tunnel data.

Table \ref{TblHarrisVars} summarizes the Harris model parameters differentiating what is known and what is calculated. The Harris PSD can be used to generate spatio-temporal wind data as described in Sec.~\ref{SecSpatialDist}.

\begin{table}
	\begin{center}
		\caption{\label{TblHarrisVars} Variable descriptions for the Harris wind-over-water turbulence model}
		\begin{tabular}{|p{0.8in}|p{2.375in}|p{1in}|p{1in}|} \hline
		{\bf{Variable}} & {\bf{Description}} & {\bf{Known/}} & {\bf{Value}} \\ 
		& & {\bf{Calculated}} & \\ \hline \hline
		$\Phi_h$ & Harris PSD (m$^2$/s) & Calculated & Eq.~\ref{EqHarrisModel}\\ \hline
		$L$ & Scaling length (m) & Known & 50 m - 126 m for altitudes 18 m - 182 m \\ \hline
		$f^*_h$ & Non-dimensional frequency & Calculated & Eq.~\ref{EqHarrisNonDimFreq} \\ \hline
		$f$ & Frequency (Hz) & Known & $f=$0-8 Hz \\ \hline 
		$C$ & Drag coefficient & Known & 0.002 or 0.0015 \\ \hline
		$u_{10}$ & Mean wind speed at 10 ft (3 m) altitude (m/s) & Known & Simulation Dependent \\ \hline
		\end{tabular}
	\end{center}
\end{table} 

\section{Forristall Model}
\label{SubSecForristall}
Forristall \cite{Forristall1988} presents an alternate spectral wind profile with data compiled from hurricanes in the Gulf of Mexico and gale force winds in the North Sea. The NWS defines gale force winds from 17.5 - 24.2 m/s (39.1 - 54.1 mph) and hurricane force winds at $>$ 33 m/s ($>$ 74 mph). Both weather events are beyond the proposed operating range for the Ascending Technologies Hummingbird UAVs; however, gusts may reach into the gale force wind range. 

Forristall concluded that the data could be fit to a Gaussian distribution defined by  
\begin{equation}
	\label{EqForristallModel}
	\Phi_f = \frac{A_f f_f^* \sigma^2}{f \left(1+B_f f^*_f\right)^{5/3}}
\end{equation}
\noindent where $A_f$ and $B_f$ are parameters determined from experimental data, $f^*_f$ is the non-dimensional frequency (different from the Harris model), $\sigma^2$ is the variance of wind speed, and $f$ is the frequency of the wind in Hz. The non-dimensional frequency and variance are defined by
\begin{align}
	\label{EqForristallNonDimFreq}
	f^*_f &= \frac{f z}{u(z)} \\
	\label{EqForristallVar}
	\sigma &= 1.92 u_*
\end{align}
\noindent where $z$ is the altitude, $u(z)$ is the wind speed at altitude $z$ as calculated by Eq.~\ref{EqLogWindProfile}, and $u_*$ is the friction velocity either calculated by Eq.~\ref{EqFrictionVelocity} or simultaneously with the sea roughness as described in Sec.~\ref{SecFrictionVelCal}. 

Table \ref{TblForristallVars} summarizes the Forristall model parameters differentiating what is known and what is calculated. The Forristall PSD can be used to generate spatio-temporal wind data as described in Sec.~\ref{SecSpatialDist}. 

\begin{table}
	\begin{center}
		\caption{\label{TblForristallVars} Variable descriptions for the Forristall wind-over-water turbulence model}
		\begin{tabular}{|p{0.8in}|p{2.375in}|p{1in}|p{1in}|} \hline
		{\bf{Variable}} & {\bf{Description}} & {\bf{Known/}} & {\bf{Value}} \\ 
		& & {\bf{Calculated}} & \\ \hline \hline
		$\Phi_f$ & Forristall PSD (m$^2$/s) & Calculated & Eq.~\ref{EqForristallModel}\\ \hline
		$A_f$ & Constant & Known & $42\pm9.07$ \cite{Forristall1988} \\ \hline
		$B_f$ & Constant & Known & $63\pm13.6$ \cite{Forristall1988} \\ \hline
		$f^*_f$ & Non-dimensional frequency & Calculated & Eq.~\ref{EqForristallNonDimFreq} \\ \hline
		$\sigma$ & Wind speed variance (m$^2$/s$^2$) & Calculated & Eq.~\ref{EqForristallVar} \\ \hline
		$f$ & Frequency (Hz) & Known & $f=$0-8 Hz \\ \hline 
		$z$ & UAV altitude (m)  & Known & UAV Dependent \\  \hline
		$u(z)$ & Wind speed at altitude $z$ (m/s) & Calculated & Eq.~\ref{EqLogWindProfile} \\ \hline
		$u_*$ & Friction velocity (m/s) & Calculated & Eq.~\ref{EqFrictionVelocity} or Sec.~\ref{SecFrictionVelCal} \\ \hline
		$u_{10}$ & Mean wind speed at 10 ft (3 m) altitude (m/s) & Known & Simulation Dependent \\ \hline
		\end{tabular}
	\end{center}
\end{table}

\section{Ochi and Shin Model}
\label{SubSecOchiShin}
The third model for comparison is the Ochi and Shin \cite{Ochi1988} model that was published in the late 1980s and is based on fitting the averaged wind data for seven different data sets. When compared to previous models, it shows much greater accuracy to the averaged data, particularly at low frequencies \cite{Chakrabarti1990}. Further supporting this, Myrhaug and Ong \cite{Myrhaug2009} show good agreement for the Ochi and Shin model compared to much more recent data taken by Andersen and L\o vseth \cite{Andersen2006} off the western coast of Norway in gale and high wind conditions. The improved accuracy comes from the piece-wise definition based on the frequency. As with the Forristall model, this model was fit to high wind data sets for use in offshore engineering structural analysis, so the accuracy of the extrapolation to lower wind speed is unknown. The PSD for the Ochi and Shin model is defined as
\begin{equation}
	\label{EqOchiShinModel}
	\Phi_{os} = \frac{\Phi_{os}^* u_*^2}{f}
\end{equation}
\noindent where $\Phi_{os}^*$ is the non-dimensional PSD and $f$ is the frequency. The non-dimensional PSD is defined by
\begin{equation}
	\label{EqOchiShinSDNonDim}
	\Phi_{os}^* = \left\{ \begin{array}{l l}
				583 f^*_f & f_f^* \leq 0.003 \\
        		\frac{420 f_f^{*^{0.7}}}{\left(1+f_f^{*^{0.35}}\right)^{11.5}} & 0.003 < f_f^* \leq 0.1 \\
            	\frac{838 f_f^*}{\left(1+f_f^{*^{0.3}}\right)^{11.5}} & f_f^* > 0.1
            	\end{array}
            	\right.
\end{equation}
\noindent where the non-dimensional frequency, $f^*_f$, is defined by Eq.~\ref{EqForristallNonDimFreq}. Table \ref{TblOchiShinVars} summarizes the Ochi and Shin model parameters differentiating what is known and what is calculated. The Ochi and Shin PSD can be used to generate spatio-temporal wind data as described in Sec.~\ref{SecSpatialDist}.

\begin{table}
	\begin{center}
		\caption{\label{TblOchiShinVars} Variable descriptions for the Ochi and Shin wind-over-water turbulence model}
		\begin{tabular}{|p{0.8in}|p{2.375in}|p{1in}|p{1in}|} \hline
		{\bf{Variable}} & {\bf{Description}} & {\bf{Known/}} & {\bf{Value}} \\ 
		& & {\bf{Calculated}} & \\ \hline \hline
		$\Phi_{os}$ & Ochi and Shin PSD (m$^2$/s) & Calculated & Eq.~\ref{EqOchiShinModel}\\ \hline
		$\Phi_{os}^*$ & Non-dimensional PSD & Calculated & Eq.~\ref{EqOchiShinSDNonDim} \\ \hline
		$f^*_f$ & Non-dimensional frequency & Calculated & Eq.~\ref{EqForristallNonDimFreq} \\ \hline
		$\sigma$ & Wind speed variance (m$^2$/s$^2$) & Calculated & Eq.~\ref{EqForristallVar} \\ \hline
		$f$ & Frequency (Hz) & Known & $f=$0-8 Hz \\ \hline 
		$z$ & UAV altitude (m)  & Known & UAV Dependent \\  \hline
		$u(z)$ & Wind speed at altitude $z$ (m/s) & Calculated & Eq.~\ref{EqLogWindProfile} \\ \hline
		$u_*$ & Friction velocity (m/s) & Calculated & Eq.~\ref{EqFrictionVelocity} or Sec.~\ref{SecFrictionVelCal} \\ \hline
		\end{tabular}
	\end{center}
\end{table}

\section{Friction Velocity Calculation}
\label{SecFrictionVelCal}
The friction velocity, $u_*$, is used in all three models described in Secs.~\ref{SubSecHarris} to \ref{SubSecOchiShin}. It is required for the drag coefficient (Eq.~\ref{EqHarrisModel}) for the Harris model, the variance (Eq.~\ref{EqForristallVar}) for the Forristall model, and the non-dimensional PSD (Eq.~\ref{EqOchiShinSDNonDim}) for the Ochi and Shin model. The calculation of the friction velocity is done in two primary ways: the Charnock model or the Volkov model \cite{Myrhaug2009}. Both calculations use the rearrangement of the logarithmic wind profile in Eq.~\ref{EqLogWindProfile} to solve for the friction velocity at an altitude of 10 m as given by
\begin{equation}
	\label{EqUStarAt10m}
	u_* = \frac{0.4 u(z)}{\ln \frac{z}{z_0}} = \frac{0.4 u_{10}}{\ln \frac{10}{z_0}}
\end{equation}

This equation is dependent on the sea roughness length, $z_0$, which is computed differently for the Charnock and Volkov models. Regardless of the model used, the sea roughness is dependent on the friction velocity, so the two variables are solved iteratively. 

The Charnock model was proposed in the mid 1950s and is only dependent on the friction velocity as follows:
\begin{equation}
	\label{EqSeaRoughnessCharnock}
	z_0 = 0.0144 u_*^2 g
\end{equation}
\noindent where $g = 9.81$ m/s$^2$.

The Charnock model is considered less accurate than the Volkov model because it does not take into account the wave age, $x_{wa} = c_p/u_*$, where $c_p$ is the phase speed of the wind waves. The wave age is ``a measure of the time the wind has been acting on a wave group" \cite{AMSWaveAge}. Since the Charnock model does not take wave age into account, Brown and Swail \cite{Brown1991} consider the Charnock model to be ``inadequate to accurately model the wind drag over water'', indicating the Volkov model should be used. 

The Volkov model was proposed in 2001 and is a piecewise function based on the wave age. For realistic conditions, the wave age reaches a maximum when the wind and waves are moving at the same rate. The model defines a non-dimensional sea roughness as follows:
\begin{equation}
	\label{EqSeaRoughnessVolkov}
	z_0^* = \left\{\begin{array}{l l}
			0.0185 & x_{wa} \leq 0.35 \\
            0.03 x_{wa} e^{-0.14 x_{wa}} & 0.35 < x_{wa} < 35 \\
            0.008 & x_{wa} \geq 35
            \end{array}
            \right.
\end{equation}
\noindent where the sea roughness length is a function of the non-dimensional sea roughness length:
\begin{equation}
	\label{EqSeaRoughnessNonDim}
	z_0 = \frac{z_0^* u_*^2}{g}
\end{equation}

To iteratively solve for the sea roughness and the friction velocity, an initial estimate of one of the parameters is needed. Brown and Swail \cite{Brown1991} present an approximation for the drag coefficient that is only dependent on a known parameter, $u_{10}$, given by  
\begin{equation}
	\label{EqCdEstimate}
	C = 1000 \left(0.65 + 0.067 u_{10}\right)
\end{equation}

The estimated drag coefficient is used to calculate an initial estimate of the friction velocity as follows:
\begin{equation}
	u_* = \sqrt{C} u_{10}
\end{equation}
To solve for the sea roughness, an estimate for the wave phase speed, $c_p$, is also needed. The realistic limits presented in \cite{Myrhaug2009} of $0.03 \leq c_p/u_{10} \leq 1.0$ are used to choose a value of $c_p$. The simulations use $c_p = 0.5 u_{10}$.

Using the initial estimate of $C$ to solve for $u_*$ and the estimate of $c_p$, the sea roughness is solved from Eq.~\ref{EqSeaRoughnessNonDim}. This sea roughness value updates the friction velocity in Eq.~\ref{EqUStarAt10m}. This process is repeated until the solutions converge. Algorithm \ref{AlgDetFrictionVel} gives the pseudocode to solve for the friction velocity, where $\varepsilon = 0.001$ for the simulations.

\begin{algorithm}
\caption{Determine the friction velocity}\label{AlgDetFrictionVel}
\begin{algorithmic}[1]
\Procedure{$[u_*,z_0,C]$ = GetFrictionVelocity($u_{10}$)}{}
\State $C \leftarrow 1000 \left(0.65 + 0.067 u_{10}\right)$
\State $u_* \leftarrow \sqrt{C} u_{10}$
\State $c_p \leftarrow 0.5 u_{10}$
\State $e_{u_*} \leftarrow 1$
\While {$e_{u_*} > \varepsilon$}
\State $z_0[i] \leftarrow \frac{z_0^* u_*[i]^2}{g}$
\State $u_*[i+1]\leftarrow\frac{0.4 u_{10}}{\ln \frac{10}{z_0[i]}}$
\State $e_{u_*} \leftarrow |u_*[i+1] - u_*[i]|$
\State $i\leftarrow i + 1$
\EndWhile
\State $C_d \leftarrow u_*^2/u_{10}^2$
\EndProcedure
\end{algorithmic}
\end{algorithm}
The Volkov model is chosen based on the results of Brown and Swail \cite{Brown1991} and Myrhaug \cite{Myrhaug2009}, who both identify the importance of including the wave age in the friction velocity calculation. In general, the Volkov model increases the PSD of the turbulence over the Charnock model \cite{Myrhaug2009}, which is more constraining for the control design. The Volkov model is used when computing the spectral model comparisons in Sec.~\ref{SubSecPSDCompare}.

\section{Spectral Model Comparison}
\label{SubSecPSDCompare}
\begin{figure}
	\begin{center}
		\includegraphics[width=6in]{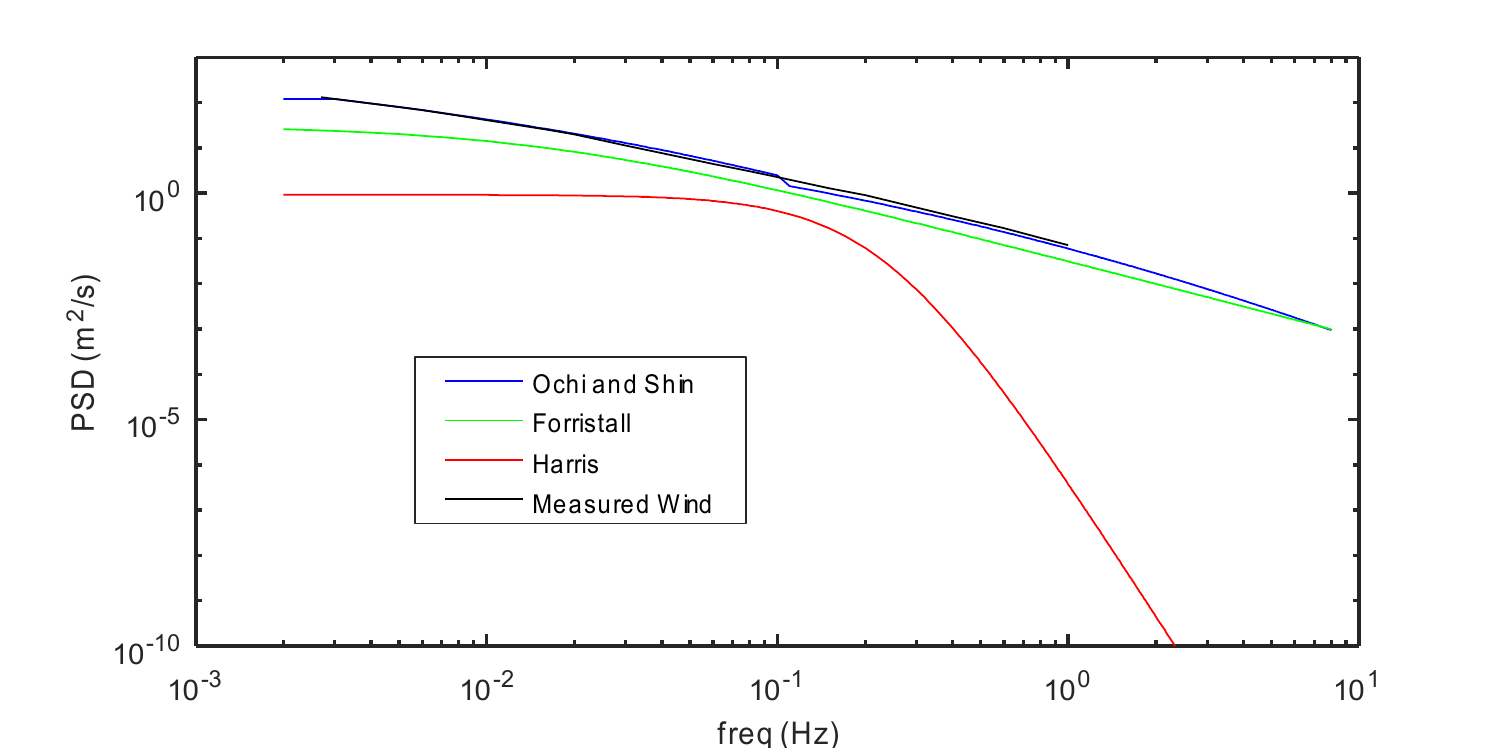}
		\caption[Comparison of the PSD for the Ochi and Shin, Forristall, and Harris models]{\label{FigPSDHFOD} Comparison of the PSD for the Ochi and Shin, Forristall, and Harris models compared to average measured wind data from \cite{Ochi1988}. As expected, the Ochi and Shin model matches the data the closest. The Forristall model is also close over a wide frequency range, but underestimates the low frequency components.}
	\end{center}
\end{figure}

The three spectral models (Harris, Forristall, and Ochi and Shin\footnote{The measured wind data from \cite{Ochi1988} is given as non-dimensional data that is converted using $\Phi_{data} = \Phi_{data}^* u_*^2/f$ and $f = f^* u(z)/z$ for $z = 10$ m, $u(z) = 10$ m/s, $c_p = 5$, and $f =0$ to $8$ Hz.}) defined in Secs.~\ref{SubSecHarris} to \ref{SubSecOchiShin} are compared to one another and the average of measured spectrum from the seven data sets in \cite{Ochi1988}. Figure \ref{FigPSDHFOD} shows the comparison. As expected, the Ochi and Shin model matches the average data very well. The Forristall model is also accurate but underestimates the low frequency PSD. The Harris model is very different from the data and other models and underestimates the PSD.

\begin{figure}
	\begin{center}
		\includegraphics[width=6in]{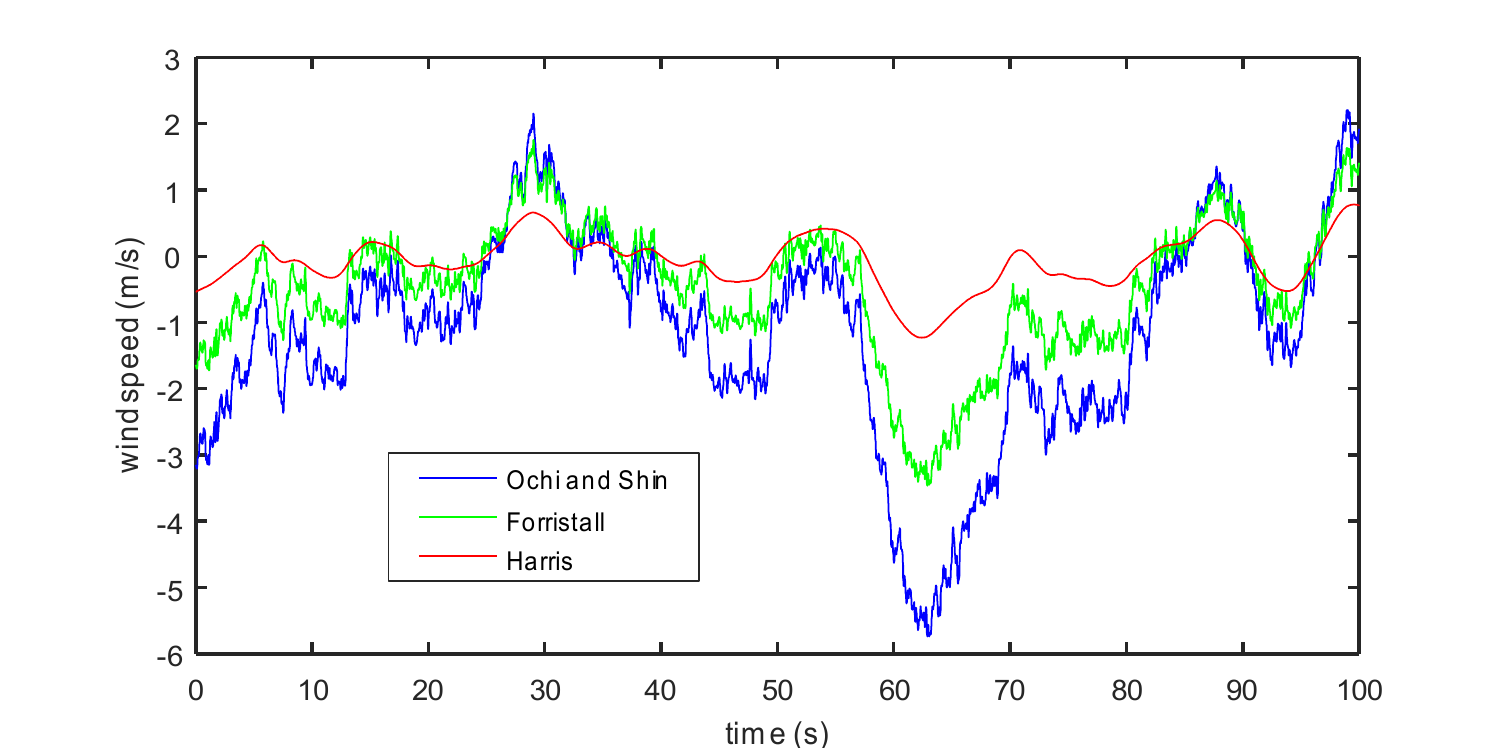}
		\caption[Comparison of time series data for the Ochi and Shin, Forristall, and Harris]{\label{FigWindHFO} Comparison of time series data for the Ochi and Shin, Forristall, and Harris models. The Ochi and Shin model shows the largest velocity magnitudes. The time series data was generated for $z = 10$ m, $u(z) = 10$ m/s, and $f =0$ to $8$ Hz. The solution to the friction velocity uses $c_p = 5$, producing $u_* = 0.45$ m/s.}
	\end{center}
\end{figure}

A comparison of the time series data for the three models is shown in Fig.~\ref{FigWindHFO}. The Ochi and Shin model shows the largest velocity magnitude, and both the Ochi and Shin and Forristall models show high frequency content for the wind model. The Harris model shows the smoothest time series because its PSD is lower at higher frequencies. The time series are generated by 
\begin{equation}
	\label{EqTimeSeries}
	y(t) = \sum_{i=1}^n \sqrt{S(\omega_i) \delta \omega} \cos \left(\omega_i t + \psi_i\right)
\end{equation}
\noindent where $S(\omega_i)$ is the PSD, $\omega$ is the frequency, and $\psi$ is a randomly generated phase angle between 0 and $2\pi$ radians \cite{Hoblit1988,Abrous2015,McMinn1997}. 

Based on the comparison of the PSD and time series data for the three models, the Ochi and Shin model is chosen for its accuracy to the averaged experimental data. It also produces the largest magnitudes for the wind speeds, which are more constraining for the control design.

\clearpage
\chapter[Spatial Wind Distribution]{5 Spatial Wind Distribution}
\label{SecSpatialDist}
The generation of a 3D spatial wavefield distribution from directional spectra is addressed by several groups \cite{Sullivan,Fossen1994,Kelly2006,Venas1998,Guo2011} due to the availability of data for directional spectra. Of the works reviewed, there were none that specifically addressed the generation of a 3D spatial windfield from spectral models; however, the principles for generating 3D wavefields can be applied with minor modifications.

In order to account for the propagation of the wind in two dimensions, a spreading function \cite{Chakrabarti1990,Forristall1997} is used. There are several spreading functions; one of the most common is the $\cos^{2s}$ function \cite{Chakrabarti1990,Guo2011} given by
\begin{equation}
	\label{EqSpreadingFunc}
	D(\omega,\theta) = D_0\cos^{2s} \theta
\end{equation}
\noindent where $s$ is a spreading parameter that is a function of frequency. The spreading function must satisfy 
\begin{equation}
	\label{EqSpreadingFuncProp}
	\int_{-\pi}^{\pi} D(\omega,\theta) d(\theta) = 1
\end{equation}
\noindent so that there is no change in the energy of the spectrum as a result of the spreading function. Typically, the spreading function is assumed to be symmetric about the mean propagation direction, where $\theta = \pm \pi/2$, and outside of this range the function is 0. Therefore, for the case $s=1$, the coefficient $D_0$ is solved using
\begin{equation}
	D_0 = \frac{1}{\int_{-\pi/2}^{\pi/2} \cos^{2s} \theta d\theta} = \frac{2}{\pi}
\end{equation}

The spreading parameter, $s$, can be specified as a constant or fit to experimental data \cite{Army1985}. While $s$ can be solved as a function of a Fourier sum of the angular spread away from the mean propagation direction from the experimental data \cite{Mitsuyasu1975,Hasselmann1980}, it is typically taken to be a constant \cite{Army1985}. In general, increased values of $s$ tend to narrow the distribution (i.e.~there is more spreading at $\theta = 0$ than $\theta = \pi/2$), so modeling of swells tends to use larger $s$. Typical values of $s$ are in the range of 1 to 10 \cite{Chakrabarti1990} or 1 to 16 \cite{Mitsuyasu1975}. For the applications considered in this dissertation, $s=1$ is chosen as this is a common value \cite{Guo2011} and does not arbitrarily narrow the propagation. 



The spreading function is multiplied by the wind spectra to form the two-dimensional spectrum as a function of frequency and spreading angle as given by
\begin{equation}
	\label{Eq2DSpectra}
	S(\omega,\theta) = S(\omega) D(\omega,\theta) = S(\omega) D(\theta) = S(\omega) D_0 \cos^2 \theta
\end{equation}
\noindent where $S(\omega)$ is the wind PSD. The superposition used to create the unidirectional time series in Eq.~\ref{EqTimeSeries} is expanded to Eq.~\ref{EqSpatialDist} where there is a double summation over the frequency range and the spreading angle. The resulting wind magnitude is
\begin{align}
	\label{EqSpatialDist}
	&v_{air,spatial}(x,y,t) = \nonumber \\
	&~~\sum_{i = 1}^{n} \sum_{j = i}^{m} \sqrt{2D_0 S(\omega_i)(\cos^2\theta_j) \Delta \omega \Delta \theta} \cos\left(k_i x \cos \theta_j + k_i y \sin \theta_j - \omega_i t + \psi_{ij}\right)
\end{align}
\noindent where $\omega_i$ varies from the lowest, $\omega_0$, to highest, $\omega_{\infty}$, frequencies that are applicable for the distribution, $\theta_j$ varies from $\pi/2$ to $\pi/2$, $k_i$ is the wavenumber, and $\psi_{ij}$ is the phase angle. The wavenumber is the spatial frequency of the wave propagation. The wavenumber for the wavefields is obtained from the linear dispersion relation for ocean waves, which is $\omega^2 = gk$ \cite{Sullivan}. The wave age calculated in Sec.~\ref{SecFrictionVelCal} shows that the wind and waves do not necessarily propagate at the same rate. The wavenumber is estimated from the frequency range of the wind (based on the literature) and the wind speed, $u_{10}$. This gives the spatial frequency as $k$ = $\omega/u_{10}$ \cite{Knight2004}. The phase angle, just as in Eq.~\ref{EqTimeSeries}, is a random value between 0 and $2\pi$. 

\begin{figure}
	\begin{center}
		\includegraphics[width=6in]{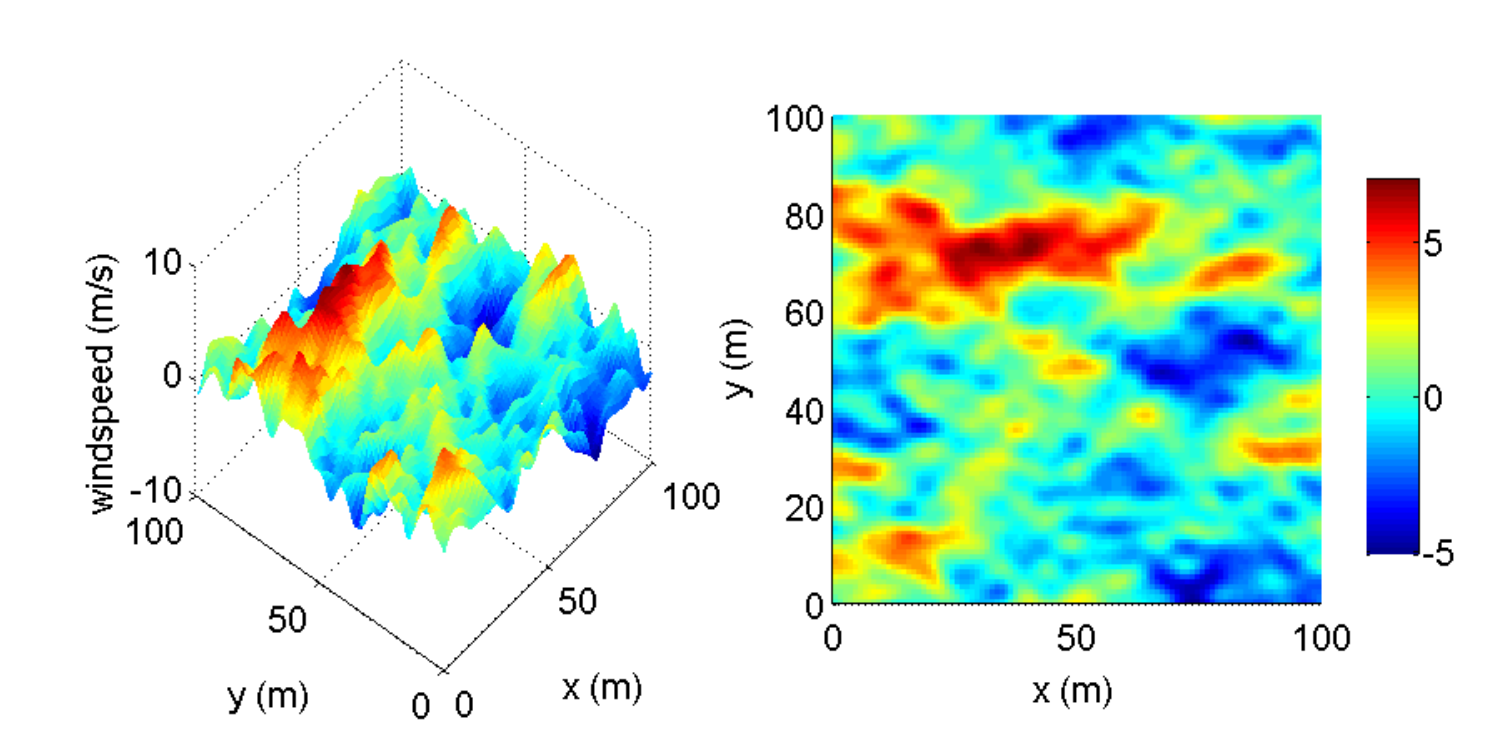}
		\caption[Spatial wavefield for input parameters $u_{10}$ = 12.8 m/s, $c_p$ = 7.5 m/s]{\label{FigOchiShinSpatial} Spatial wavefield for input parameters $u_{10}$ = 12.8 m/s, $c_p$ = 7.5 m/s, at an altitude at $z$ = 2.5 m.}
	\end{center}
\end{figure}

\begin{figure}[h!]
	\begin{center}
		\includegraphics[width=6in]{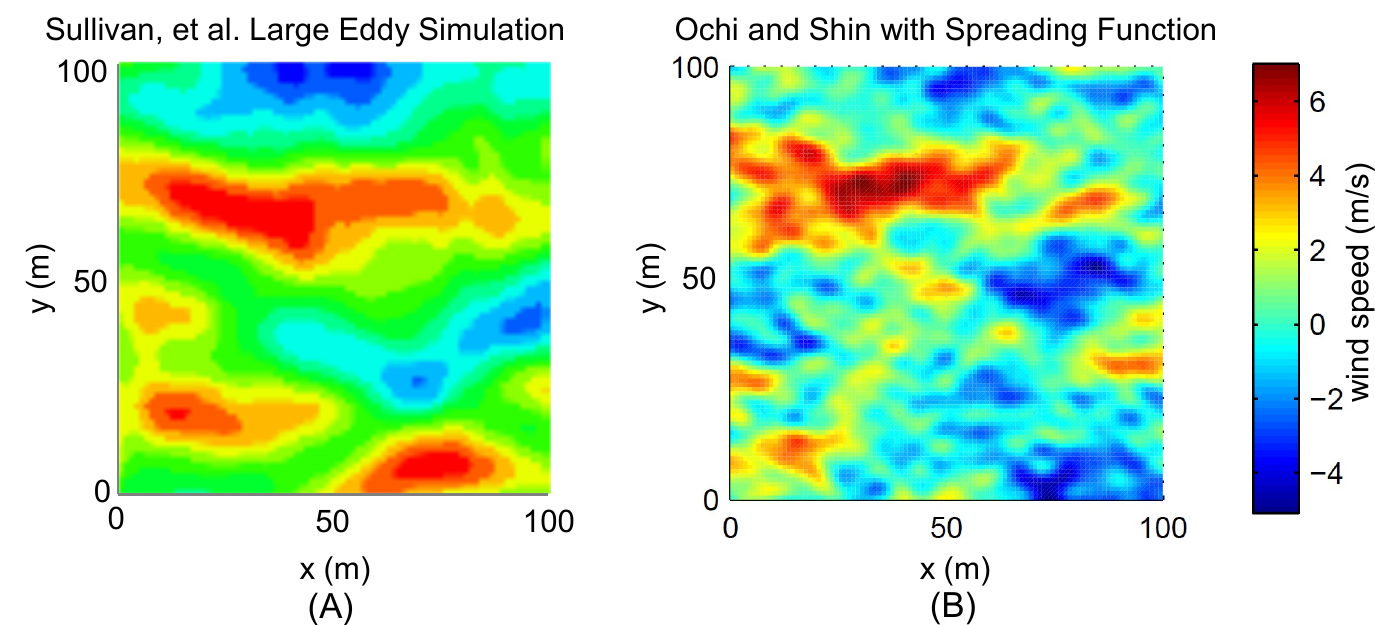}
		\caption[Comparison of Sullivan and Ochi and Shin spectral distributions]{\label{FigCompareSullivan1vsOchi} Comparison of (A) Sullivan (Figure 5 right panel from \cite{Sullivan}) (B) Ochi and Shin spectral spatial distributions over a 100 m x 100 m area for starting parameters $u_{10}=12.8$ m/s, $z = 2.5$ m, and $c_p/u_* \approx 1.4$ m/s.}
	\end{center}
\end{figure}

An example of the spatial distribution generated by Eq.~\ref{EqSpatialDist} is shown in Fig.~\ref{FigOchiShinSpatial} using the Ochi and Shin model. A qualitative comparison of the spreading function wind field generation and the large eddy simulation high fidelity model is shown in Fig.~\ref{FigCompareSullivan1vsOchi} for the same set of starting parameters. The relative periodicity and spatial size is fairly well re-created, confirming Eq.~\ref{EqSpatialDist} is a reasonable method for spatio-temporal wind field generation.  

\section{Wind Field Propagation}
\begin{figure}
	\begin{center}
		\includegraphics[width=3in]{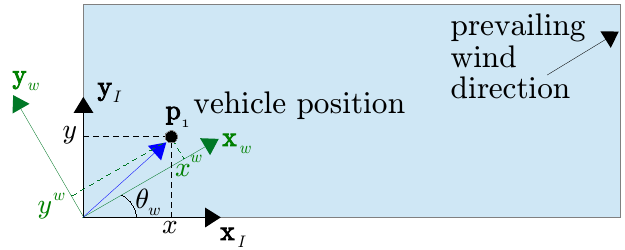}
		\caption[Coordinate transformation to simulate wind fields in any direction.]{\label{FigRotateWind} Coordinate transformation to simulate wind fields in any direction. The prevailing wind is in the $\mathbf{x}_w$ direction, and the wind magnitude is computed for any point $\mathbf{p}_1$.}
	\end{center}
\end{figure}

The definition in Eq.~\ref{EqSpatialDist} generates a spatial wind field that propagates with time in the $+\mathbf{x}_I$ direction, as shown in Fig.~\ref{FigRotateWind}. If the prevailing wind is in a different direction, the field must be rotated using a coordinate transformation. The new points in the frame, $\mathbf{x}_w$ and $\mathbf{y}_w$, are defined as $\mathbf{p}_1^w = [x^w,~y^w,~z^w]$. This is illustrated in Fig.~\ref{FigRotateWind} and used to determine the turbulence magnitude by: 
\begin{align}
	\label{EqSpatialDistCoordTransform}
	v_{turb}(x,y,t) &= \sum_{i = 1}^{n} \sum_{j = 1}^{m} \sqrt{2 S(\omega_i,\theta_j)\Delta \omega \Delta \theta} \cos\left(k_i x^w \cos \theta_j + k_i y^w \sin \theta_j - \omega_i t + \psi_{ij}\right)
\end{align}

This rotation is important since MIL-F-8785C states that for the low altitude model, the longitudinal component shall be aligned with the prevailing wind direction. 

\section{Vector Components}
\label{SubSecVectorComponents}
In addition to the rotation of the wind field to align with the prevailing wind direction, the scalar magnitude produced by Eq.~\ref{EqSpatialDist} must be translated into a wind vector. The Von K\'arm\'an model has four PSD components that are all independent; however, the Ochi and Shin model has only a single PSD function. 

In the Von K\'arm\'an model the PSD functions from Tbl.~\ref{TblVKLAM} are taken to be the independent vector components. Since the longitudinal component is aligned with the prevailing wind direction, this defines the relationship between the longitudinal and lateral components in the inertial frame, as shown in Fig.~\ref{FigRotateWind}. Equations \ref{EqVKTransVector} and \ref{EqVKRotVector} summarize the vector implementation of the full wind velocity experienced by the vehicle.
\begin{align}
	\label{EqVKTransVector}
	\mathbf{v}_{air} &= \left[\begin{array}{c}
	\left(v_{turb,lon}(x,y,t) + v_{air} + v_{gust}\right)\cos(\theta_w) + v_{turb,lat}(x,y,t)\sin(\theta_w)  \\
	\left(v_{turb,lon}(x,y,t) + v_{air} + v_{gust}\right)\sin(\theta_w) + v_{turb,lat}(x,y,t)\cos(\theta_w)   \\
	v_{turb,vert}(x,y,t)
	\end{array}
	\right]					
\\
	\label{EqVKRotVector}
	\boldsymbol\omega_{air}& = \left[\begin{array}{c}
	\omega_{air,roll} \\
	\omega_{air,pitch} \\
	\omega_{air,yaw}
	\end{array}
	\right]
\end{align}
\noindent where $v_{turb,lon}$, $v_{turb,lat}$, $v_{turb,vert}$, $\omega_{air,roll}$, $\omega_{air,pitch}$, $\omega_{air,yaw}$ are computed by Eq.~\ref{EqSpatialDistCoordTransform} using the respective Von K\'arm\'an PSD functions in Tbl.~\ref{TblVKLAM}, $v_{air,mean}$ is the mean wind speed, and $v_{gust}$ is the magnitude of the gust (as defined in Sec.~\ref{SecGustModel}). 

The Ochi and Shin model is not quite as straightforward since there is only a single PSD. Hsu and Blanchard \cite{Hsu2004} provide an initial starting point for scaling the turbulence intensities for the lateral and vertical components, but the Von K\'arm\'an components provide insight into the relative magnitude ratios as well. Assuming that the Ochi and Shin model is the analog of the longitudinal component of the Von K\'arm\'an model, then the ratio of the other 5 components in the Von K\'arm\'an model to the longitudinal component are used to scale the Ochi and Shin model PSD. In addition to magnitude, the correlation between the translational and rotational spectra is also considered. Again using the Von K\'arm\'an model, there are four independent components (longitudinal, lateral, vertical, roll) and two dependent components (pitch/vertical and yaw/lateral). This correlation is implemented by using the same random sequence of values for $\psi$ in Eq.~\ref{EqSpatialDistCoordTransform} for the components that are correlated. 

\begin{figure}
	\begin{center}
		\includegraphics[width=6in]{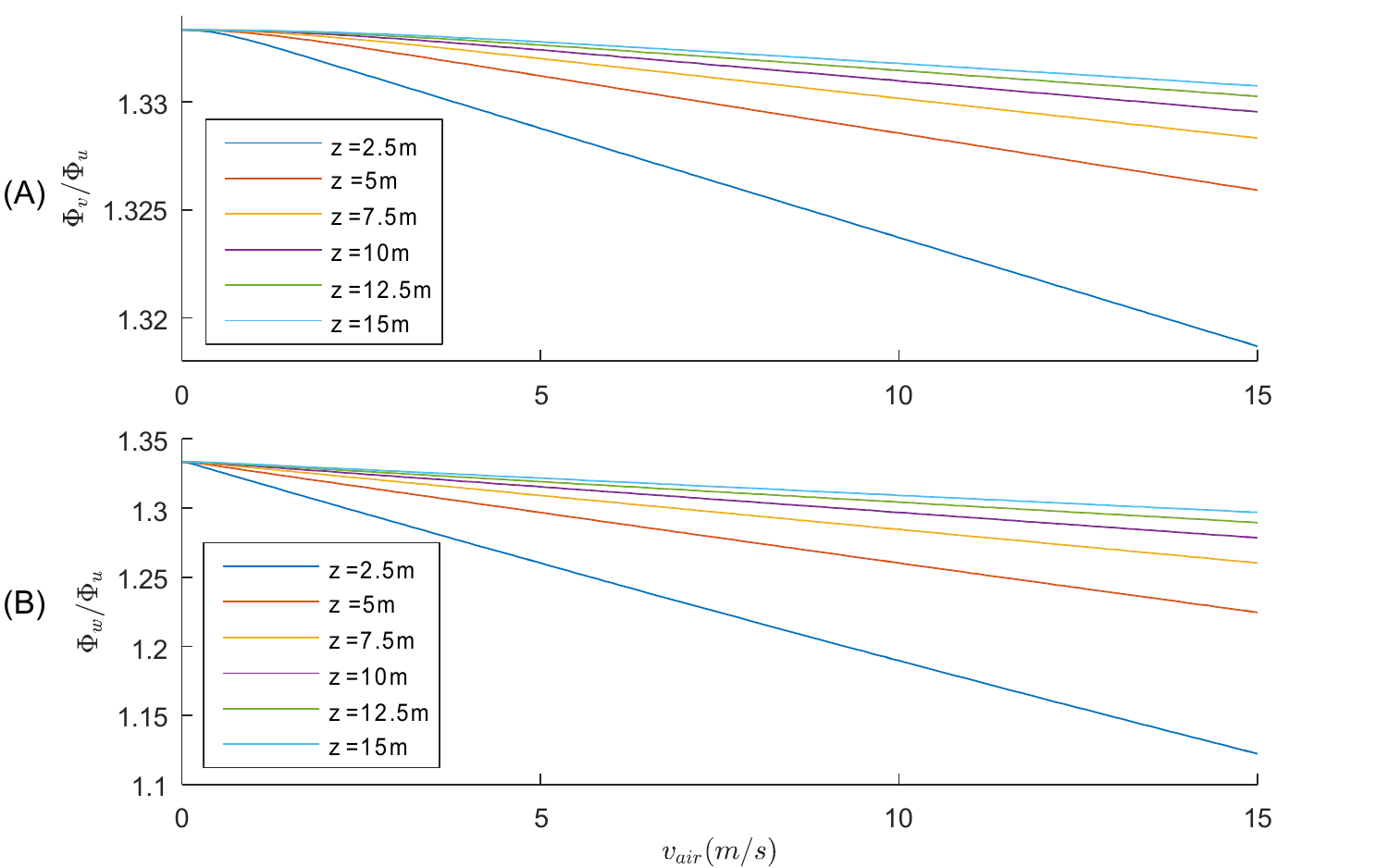}
		\caption[Comparison of the ratio of the Von K\'arm\'an translational components]{\label{FigOSRatioTrans} Comparison of the ratio of the Von K\'arm\'an (A) lateral PSD ($\Phi_{v}$) and (B) vertical PSD ($\Phi_w$) to the longitudinal PSD ($\Phi_u$) for varying mean wind speeds, $v_{air}$, and altitudes, $z$.}
	\end{center}
\end{figure}

\begin{figure}
	\begin{center}
		\includegraphics[width=6in]{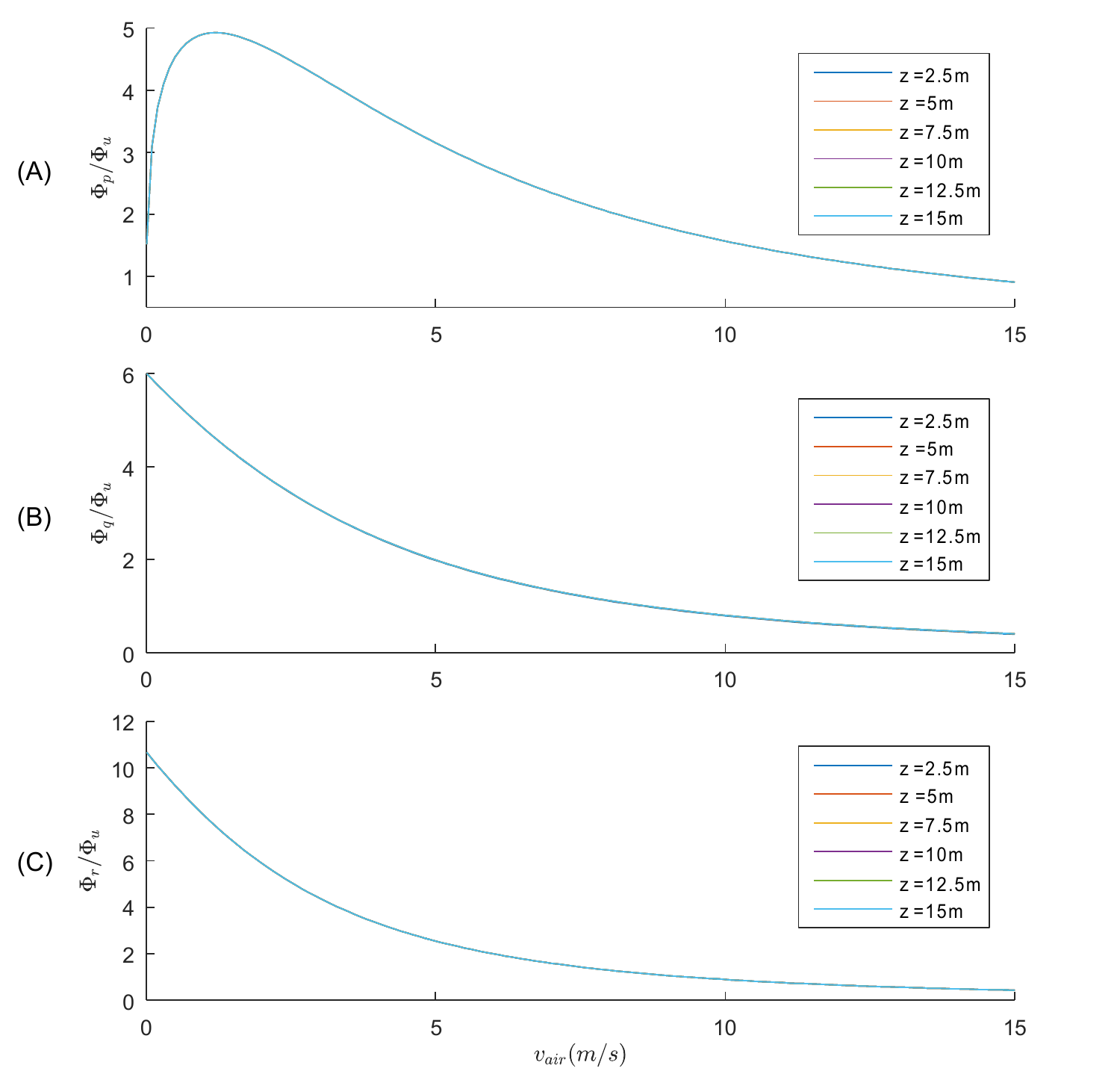}
		\caption[Comparison of the ratio of the Von K\'arm\'an rotational components]{\label{FigOSRatio} Comparison of the ratio of the Von K\'arm\'an (A) roll PSD ($\Phi_{p}$), (B) pitch PSD ($\Phi_q$), and (C) yaw PSD ($\Phi_r$) to the longitudinal PSD ($\Phi_u$) for varying mean wind speeds, $v_{air}$, and altitudes, $z$. Note that dependence on altitude is very minimal.}
	\end{center}
\end{figure}

Figures \ref{FigOSRatioTrans} to \ref{FigOSRatio} show the relative magnitudes for the translational and rotational components, respectively, for various wind speeds. Notice that the rotational components vary minimally with the altitude. These relationships are used to define the components of the Ochi and Shin model for use in Eqs.~\ref{EqVKTransVector} and \ref{EqVKRotVector}.

\chapter[Gust Model]{6 Gust Model}
\label{SecGustModel}
The final component to the spatio-temporal wind field is the gust model. Gusts are only considered for the translational components of wind as, to the best of the author's knowledge, no rotational gust models have been proposed to date. The gusting is based on wind farm statistical gust data from \cite{Branlard2009,Manwell2003}, which is limited to measuring only 2D gusts and does not account for boundary layer effects near the ground. These studies suggest the following five parameters when defining a gust:
\begin{singlespace}
\begin{enumerate}
	\item Amplitude - maximum amplitude of the gust
	\item Shape - amplitude vs time profile of the gust
	\item Duration - time span from gust start to finish
	\item Propagation speed - speed at which the gust moves in the prevailing wind direction
	\item Occurrence frequency - time between discrete gust events
	\item Longitudinal and Lateral spread/dissipation - spatio-temporal modeling of the gust
\end{enumerate}
\end{singlespace}
Each of these parameters is considered for discrete gust events and is discussed briefly in the following subsections. The combination of these parameters generates the final gust amplitude to use in Eq.~\ref{EqVKTransVector}.

\section{Amplitude}
\label{SubSecGustAmp}
\begin{figure}
	\begin{center}
		\includegraphics[width=3in]{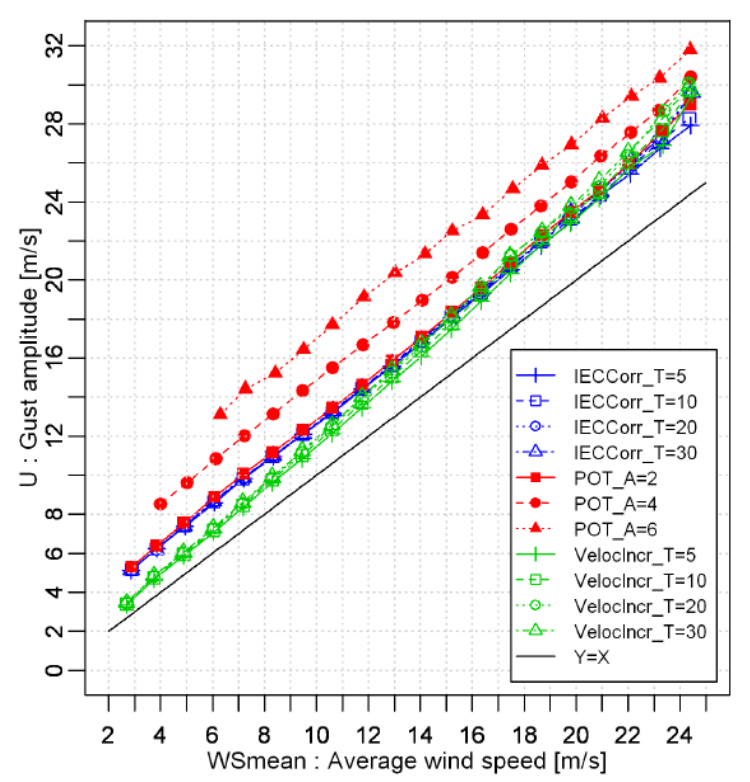}
		\caption[Relationship between the mean speed and the maximum gust amplitude]{\label{FigBranlard43} Relationship between the mean speed and the maximum gust amplitude. The relationship is nearly linear for all the different gust detection methods \cite{Branlard2009} (Figure 4.3).}
	\end{center}
\end{figure}

The gust amplitude is strongly correlated to the mean wind speed \cite{Branlard2009}. The plot in Fig.~\ref{FigBranlard43} is taken from \cite{Branlard2009} to show the strong correlation between mean wind speed and gust amplitude. A linear fit is used to define the maximum gust amplitude as a function of mean wind speed for the ``POT\_A=2" line from Fig.~\ref{FigBranlard43}. This line uses the peak over threshold (POT) method, which identifies a gust if the wind magnitude increases by the threshold amount (in this case 2 m/s for POT\_A=2) over the mean wind speed. For the vehicles considered in this research, 2 m/s changes could be significant. Additionally, the ``POT\_A=2" line fits closely with the ``IECCorr" lines, which use the correlation method that identifies gusts by matching gust shape, and the ``VelocIncr" lines, which measure the acceleration from the first to last point of a moving window to determine if a gust is present. The resulting linear fit is given by
\begin{equation}
	\label{EqGustAmplitude}
	v_{gust} = \frac{52}{43} v_{air} - \frac{1}{8}
\end{equation}
 
\section{Shape}
\label{SubSecGustShape}
\begin{figure}
	\begin{center}
		\includegraphics[width=3in]{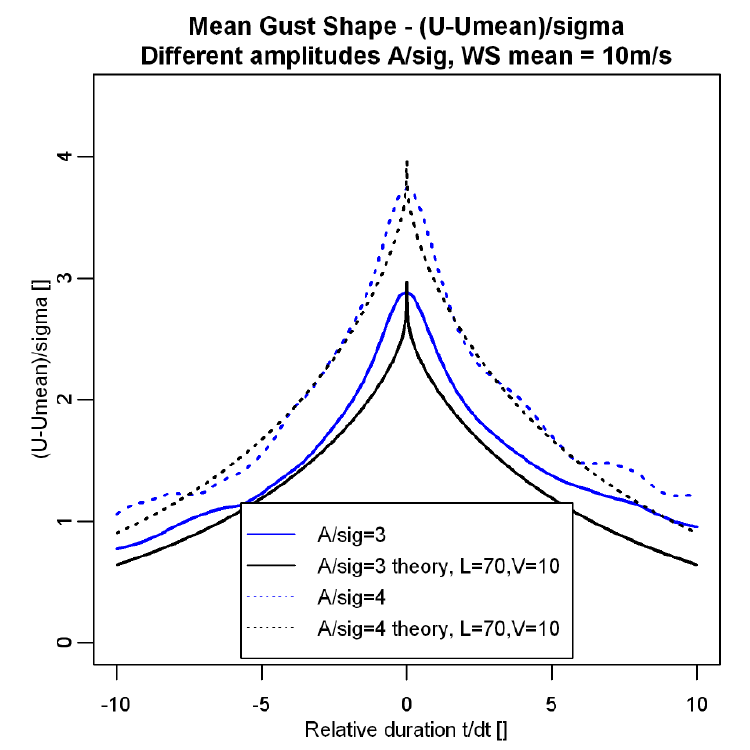}
		\caption[Mean gust shape characterization by comparing non-dimensional parameters]{\label{FigBranlard412} Mean gust shape characterization by comparing non-dimensional parameters of relative duration (data before and after the gust peak, in this case $\pm10$ seconds with $dt = 1$) and change in wind speed over standard deviation \cite{Branlard2009} (Figure 4.12). While there is good agreement of experiment and theory, this does not provide insight for generating time series gust data.}
	\end{center}
\end{figure}

It is difficult to define a generic gust shape based on the data that was collected at the wind farm \cite{Branlard2009}. No clear spatio-temporal gust shape emerges from the data. Instead, the authors characterize the trends of the gust shapes relative to the gust peak by examining the data before and after the gust peak. The resulting plot is in terms of the relative duration (in this case data $\pm10$ seconds around the peak with $dt = 1$), compared to the non-dimensional parameter, $(v_{gust}-v_{air})/\sigma_{air}$, where $\sigma_{air}$ is the standard deviation over a 10 minute period that is also used to calculate the mean wind speed, $v_{air}$. Figure \ref{FigBranlard412} shows what the authors and other works \cite{Larsen2003} call the ``mean gust shape". While the experimental data and theory agree, there is not an obvious way to generate a gust shape in the time domain from this data. 

\begin{figure}
	\begin{center}
		\includegraphics[width=3in]{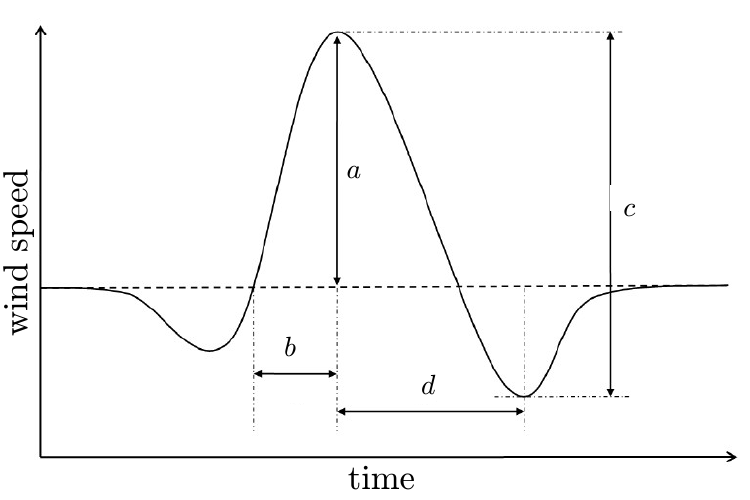}
		\caption[Representative gust shape defined by \cite{Branlard2009,Manwell2003}]{\label{FigGenericGust} Representative gust shape defined by \cite{Branlard2009,Manwell2003} provides the basis for the generic gust definition used for this dissertation.}
	\end{center}
\end{figure}

For this research, the gust shape considered is a discrete event that produces a temporary increase in wind speed over the mean that then returns to its original value. Consequently, gusts that change the overall mean wind speed are considered separately and use the $1-\cos$ shape. A generic shape for the discrete gust event, such as that shown in Fig.~\ref{FigGenericGust}, was suggested by \cite{Manwell2003} and is used as the basic gust shape for this dissertation.  

\begin{figure}
	\begin{center}
		\includegraphics[width=6in]{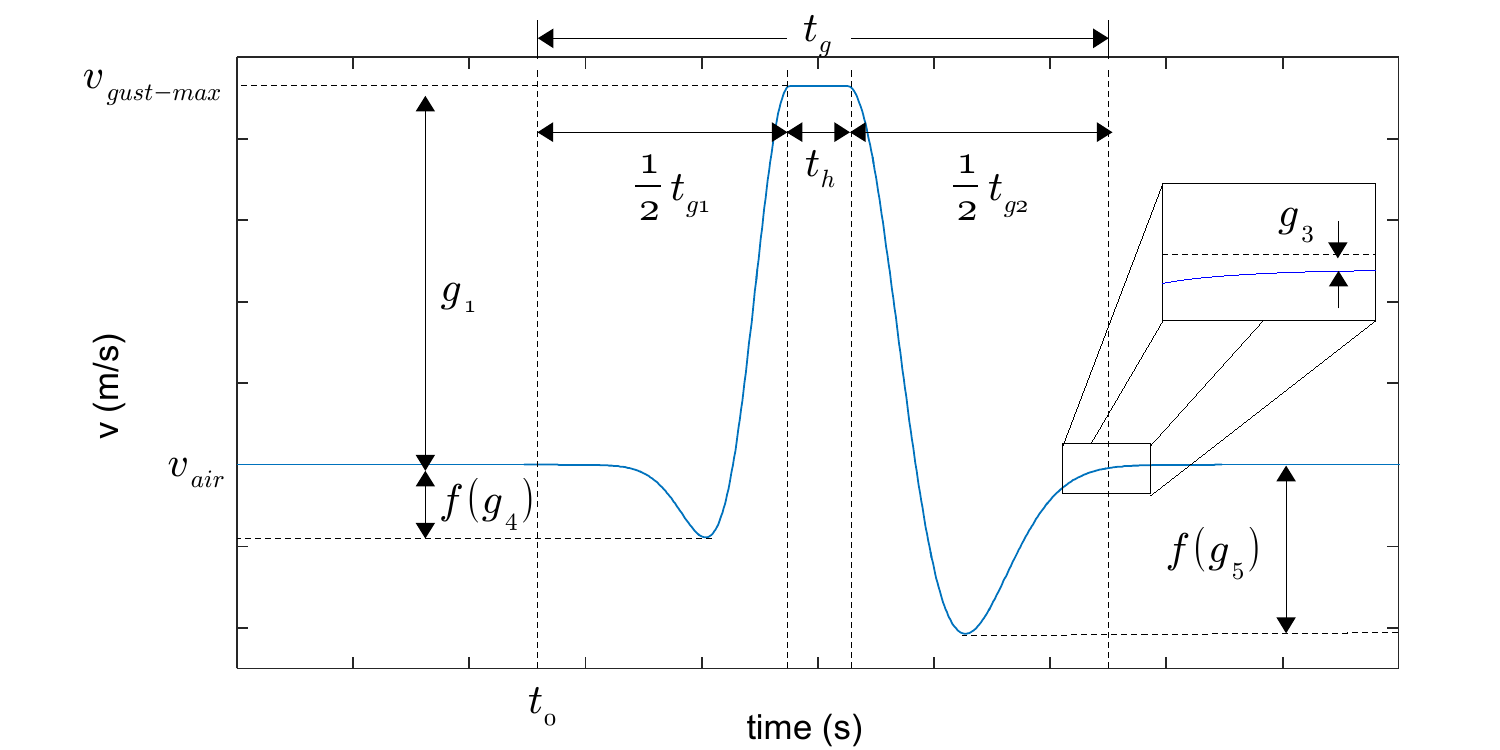}
		\caption[Representative gust shape allowing independent specification]{\label{FigGenericGustFlatTop} Representative gust shape allowing independent specification for the rising, falling, and peak hold portions of the gust.}
	\end{center}
\end{figure}

Branlard \cite{Branlard2009} provides the equation for a symmetric curve of similar shape to that in Fig.~\ref{FigGenericGust}. To allow greater control, the definition from \cite{Branlard2009} is generalized as shown in Fig.~\ref{FigGenericGustFlatTop} and defined as follows:
\begin{equation}
	\label{EqGenericGust}
	v_{gust,nom}(t) = \left\{\begin{array}{ll}
	g_1 \left(1-\left(g_{2r} t - g_3\right)^2\right) e^{-\left(g_{2r} t - g_3\right)^2/g_4}, & 0 < t < \frac{1}{2} t_{g1} \\
	g_1, & \frac{1}{2} t_{g1} \leq t \leq \frac{1}{2}t_{g1} + t_h \\
	g_1 \left(1-\left(g_{2f} t^* - g_3\right)^2\right) e^{-\left(g_{2f} t^* - g_3\right)^2/g_5}, & t > \frac{1}{2} t_{g1}  + t_h
	\end{array}
	\right.
\end{equation}
\noindent where $g_1$ is the gust amplitude obtained from Eq.~\ref{EqGustAmplitude}, $g_{2r} = 2 g_3/t_{g1}$ and $g_{2f} = 2 g_3/t_{g2}$, which scale the time for the rising and falling sides of the curve respectively, $g_3=6$ is the time, $t$, where $(1-t^2)e^{-t^2/2}$ is within 1\% (or less) of zero, $g_4 > 0$ impacts the magnitude of the ``dip" before the gust, $g_5 > 0$ impacts the magnitude of the ``dip" after the gust, $t_h$ is the time duration at the peak gust magnitude, $t^* = t + \frac{1}{2} (t_{g2}-t_{g1}) - t_h$, $t_{g1}$ is the time from gust start to peak for the rising side, and $t_{g2}$ is the time from peak to gust end for the falling side. 

\section{Duration}
\label{SubSecGustDuration}
\begin{figure}
	\begin{center}
		\includegraphics[width=3in]{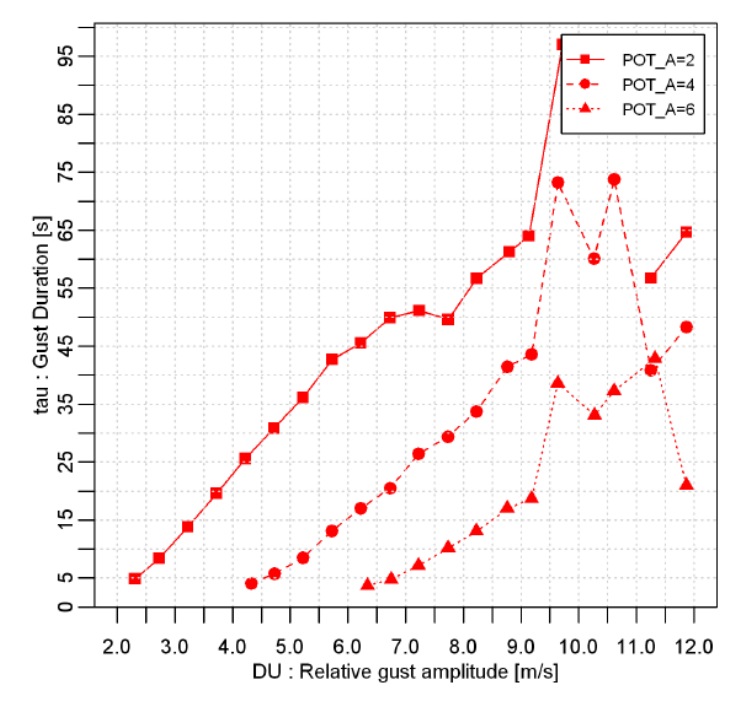}
		\caption[Relationship between relative gust amplitude and gust duration]{\label{FigBranlard410} Relationship between relative gust amplitude and gust duration where the relative gust amplitude is defined as $v_{gust} - v_{air}$. Reasonable wind speeds for this dissertation generally correspond to the linear region of the plot \cite{Branlard2009} (Figure 4.10).}
	\end{center}
\end{figure}

There is a fair amount of variation in the gust duration as presented in \cite{Branlard2009} and shown in Fig.~\ref{FigBranlard410}. The linear region of the plot is used since precise accuracy of the gust duration is not required, and the relative gust magnitude in the linear region is consistent with the wind speeds examined for this research. To be consistent with the gust amplitude linear fit, the ``POT\_A=2" line is used again. The linear fit is given by
\begin{equation}
	\tau_{gust} = 0.71\left(v_{gust} - v_{air}\right) + 3.51
\end{equation}

\section{Propagation}
The data from \cite{Branlard2009} indicates a reasonable amount of variation in the data for correlating gust propagation. There are a two parameters that stand out for consideration: (1) the presence of the expected gust (i.e. did the gust actually propagate as expected, and, if so, how did the shape of the gust compare to what was experienced at the first measurement), and (2) the propagation speed. 

There is not a good model that is gleaned from this data to realistically modify the gust propagation so that either the shape changes or the gust dissipates out. As such, the only modifications to the gust in the prevailing wind direction come from the changing turbulence as calculated by either the Von K\'arm\'an or Ochi and Shin turbulence model. Given that the formation distributions are not spanning several hundred meters or more, this assumption is reasonable (discussed further in Sec.~\ref{SubSecLatLongDissipation}). 

The propagation speed is the mean wind speed, $v_{air}$, as supported by \cite{Branlard2009}. It is assumed that the propagation speed is constant over the simulation, which is reasonable for simulations on the order of tens of minutes as the mean wind speed is typically computed over 10 minute intervals for wind data analysis, and the 10 minute averages change slowly \cite{Branlard2009}. 



\section{Occurrence Frequency}
Branlard \cite{Branlard2009} presents a summary of the number of gusts per month as detected by the various algorithms to determine how frequently gusts occur. Again considering the peak over threshold model for a rise of 2 m/s, there are 9,777 gusts detected per month. This extrapolates to 13.5 gusts per hour or a gust roughly every 4 minutes. This estimate is used as an average value for generating gusts in the simulations.  

\section{Lateral and Longitudinal Dissipation}
\label{SubSecLatLongDissipation}
\begin{figure}
	\begin{center}
		\includegraphics[width=3in]{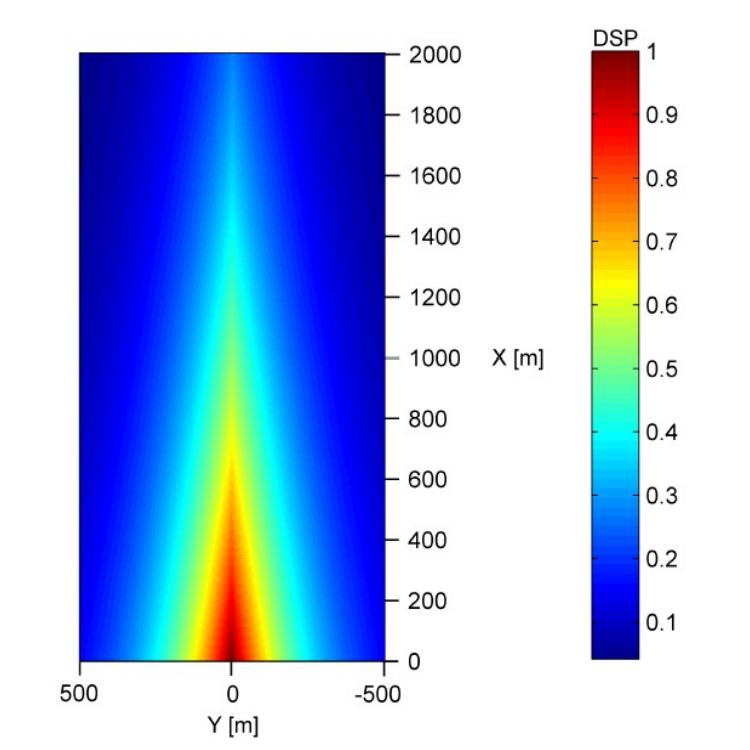}
		\caption[Probability of detecting a gust relative to a wind turbine]{\label{FigBranlard69} Probability of detecting a gust relative to a wind turbine that just detected a gust \cite{Branlard2009} (Figure 6.9).}
	\end{center}
\end{figure}

To determine the lateral and longitudinal spread of the gust, the authors of \cite{Branlard2009} go through the data manually to determine how individual gusts propagate. In general, there is very strong correlation ($\approx$85\%) for the longitudinal presence of a gust at a point in space 200 m away from but in line with the initial detection location, as shown in Fig.~\ref{FigBranlard69}. This is well within the assumed operating area so it is reasonable to seek a model that shows strong correlation longitudinally.

The lateral spread of the gust drops off significantly around points 250 m away from and perpendicular to the initial detection location based on the wind farm data as shown in Fig.~\ref{FigBranlard69}. Since the assumed extents of the formations are less than 250 m in any direction, the model should show some lateral spread of gusts to other vehicles. 

The assessment of the propagation is somewhat dependent on the locations of the wind turbines. The probabilities and correlation functions are given at specific distances and then curve fit for all distances. Additionally, since the data is measured on wind turbine towers, there is a minimum distance required to avoid the wake stream of one turbine affecting the measurement of another. The proposed curve fit \cite{Branlard2009} for probability of detecting a gust laterally and longitudinally is
\begin{equation}
	\label{EqGustSpreadFit}
	P = e^{-\frac{\delta x}{l_x}} e^{-\frac{|\delta y|}{l_y}}
\end{equation}
\noindent where $l_x$ and $l_y$ are characteristic length parameters associated with the gust and $\delta x$ and $\delta y$ are the relative distance between the ``gust center" and the detection point. The characteristic length is calculated by fitting Eq.~\ref{EqGustSpreadFit} to the data. 

Equation \ref{EqGustSpreadFit} also provides the basis for scaling the spreading and dissipation of the gust, where the probability is used as the gust magnitude scaling parameter, $K_{gust,lat} = P$. The characteristic length parameters are based on the expected wind environment and $\delta x$ and $\delta y$ are the distance between the gust center and the vehicle location along the $\mathbf{x}_w$ and $\mathbf{y}_w$ axes, respectively. The gust center propagates at a speed of $v_{air}$ in the $\theta_w$ direction, and the resulting final gust magnitude is
\begin{equation}
	\label{EqGustFinal}
	v_{gust} = K_{gust,lat} v_{gust,nom}(t^*)
\end{equation}
\noindent where
\begin{equation}
	t^* = t -  \frac{1}{v_{air}} \left((\mathbf{p}(t) - \mathbf{p}_{gust}) \cdot [\cos \theta_w,~\sin \theta_w,~0]^T\right)
\end{equation}
\noindent $\mathbf{p}_{gust}$ is the ``starting" position of the gust and $t$ is the simulation time. This gust magnitude in Eq.~\ref{EqGustFinal} is used in Eq.~\ref{EqVKTransVector} to fully define the wind model at any point $\mathbf{p}$.

\chapter[Conclusions]{7 Conclusions}
The spatio-temporal wind model defined in this chapter addresses a deficiency in the literature for realistic spatio-temporal wind field modeling for UAV formations that includes turbulence and gusting. The model meets the requirements of preserving realistic features of the environmental turbulence and gusting as well as computational feasibility. The Von K\'arm\'an model is chosen for operation in over-land environments where it has been shown that this model is realistic even for smaller UAVs. The Ochi and Shin model is one of the most accurate wind-over-water turbulence models and is chosen for the over-water scenarios. To create the spatio-temporal distribution a spreading function uses the chosen PSD functions to generate a correlated wind field where each vehicle experiences unique conditions. Lastly, the spreading function definition is also suitable for superposing with a 2D gust model. 

This environmental wind model is intended for use in simulations of groups of vehicles to simulate the unique but correlated wind conditions experienced by the vehicles. The stochastic basis of the model allows multiple wind fields to be used for the verifying control laws.

\begin{singlespace}
\addcontentsline{toc}{chapter}{References}
\bibliographystyle{unsrt}
\bibliography{C:/Users/cole/Documents/ColeThesisWork/Documentation/ReferencesBibFiles/Controls,C:/Users/cole/Documents/ColeThesisWork/Documentation/ReferencesBibFiles/EstimationFormation,C:/Users/cole/Documents/ColeThesisWork/Documentation/ReferencesBibFiles/Gusting,C:/Users/cole/Documents/ColeThesisWork/Documentation/ReferencesBibFiles/RCSBib,C:/Users/cole/Documents/ColeThesisWork/Documentation/ReferencesBibFiles/WindWaterRefs}
\end{singlespace}

\end{document}